\documentclass[sn-mathphys,Numbered]{sn-jnl}

\usepackage{mathtools}%
\usepackage{bm}%
\usepackage{natbib}%
\usepackage{bibentry}%

\usepackage{graphicx}%
\usepackage{multirow}%
\usepackage{amsmath,amssymb,amsfonts}%
\usepackage{amsthm}%
\usepackage{mathrsfs}%
\usepackage[title]{appendix}%
\usepackage{xcolor}%
\usepackage{textcomp}%
\usepackage{manyfoot}%
\usepackage{booktabs}%
\usepackage{algorithm}%
\usepackage{algorithmicx}%
\usepackage{algpseudocode}%
\usepackage{listings}%
\usepackage{todonotes}%
\usepackage{longtable}%

\theoremstyle{thmstyleone}%
\newtheorem{theorem}{Theorem}
\theoremstyle{thmstyletwo}%

\theoremstyle{thmstylethree}%

\newcommand{\bracket}[2]{[#1,#2]}

\makeatletter
\newcommand\smallscript{\@setfontsize\smallscript{7.5}{8.5}}
\makeatother

\newcommand{\firstrevision}[1]{#1}
\newcommand{\firstrevisioncolor}{}

\begin{document}

\title[Efficient Separation of RLT Cuts for Implicit and Explicit Bilinear \firstrevision{Terms}]{Efficient Separation of RLT Cuts for Implicit and Explicit Bilinear \firstrevision{Terms}\footnote{This is the extended version of the article published as: Bestuzheva, K., Gleixner, A., Achterberg, T.: Efficient separation of RLT cuts for
implicit and explicit bilinear products. In: Integer Programming and Combinatorial Optimization: 24th International Conference, IPCO 2023, Madison, WI,
USA, June 21–23, 2023, Proceedings, pp. 14–28. Springer, Berlin, Heidelberg
(2023). \url{https://doi.org/10.1007/978-3-031-32726-1_2}.}}


\author*[1]{\fnm{Ksenia} \sur{Bestuzheva}}\email{bestuzheva@zib.de}

\author[1,2]{\fnm{Ambros} \sur{Gleixner}}\email{gleixner@htw-berlin.de}

\author[3]{\fnm{Tobias} \sur{Achterberg}}\email{achterberg@gurobi.com}

\affil*[1]{\orgdiv{AI in Society, Science, and Technology}, \orgname{Zuse Institute Berlin}, \orgaddress{\street{Takustra{\ss}e 7}, \city{Berlin}, \postcode{14195}, \country{Germany}}}

\affil[2]{\orgname{HTW Berlin}, \orgaddress{\city{Berlin}, \postcode{10313}, \country{Germany}}}

\affil[3]{\orgname{Gurobi GmbH}, \orgaddress{\street{Sandstra{\ss}e 104}, \city{Monheim am Rhein}, \postcode{40789}, \country{Germany}}}


\abstract{The reformulation-linearization technique (RLT) is a prominent approach to constructing tight
linear relaxations of non-convex continuous and mixed-integer optimization problems.
The
goal of this paper is to extend the applicability and improve the performance of RLT for bilinear
product relations.
First, we present a method for detecting bilinear product relations implicitly contained in mixed-integer
linear programs, which is based on analyzing linear constraints with binary variables,
thus enabling the application of bilinear RLT to a new class of problems.
Strategies for filtering product relations are discussed and tested.
Our second contribution addresses the high computational cost of RLT cut separation, which presents
one of the major difficulties in applying RLT efficiently in practice.
We propose a new RLT cutting plane separation algorithm which identifies combinations of
linear constraints and bound factors that are expected to yield an inequality that is violated by
the current relaxation solution.
\firstrevision{This algorithm is applicable to RLT cuts generated for all types of bilinear terms, including but not
limited to the detected implicit products.}
A detailed computational study based on independent implementations in two solvers evaluates the performance impact of the
proposed methods.}




\maketitle

\section{Introduction}

The reformulation-linearization technique (RLT) was first proposed by Adams and
Sherali~\cite{adams1986tight,adams1990linearization,adams1993mixed} for linearly constrained
bilinear problems with binary variables, and has been applied to
mixed-integer~\cite{sherali1990hierarchy,sherali1994hierarchy}, semi-infinite~\cite{sherali2009reformulation},
general bilinear~\cite{sherali1992new} and polynomial~\cite{sherali1992global} problems.
RLT constructs valid polynomial constraints, then linearizes these constraints by
using polynomial relations given in the problem and applying relaxations when such relations are not
available.
If relations used in the linearization step are violated by a relaxation solution, this procedure
may yield violated cuts.
By increasing the degree of derived polynomial constraints, hierarchies of relaxations can be
constructed, which were shown to converge to the convex hull representation of \firstrevision{mixed-integer linear programming problems (MILPs)}
and mixed-integer polynomial problems where continuous variables appear
linearly~\cite{sherali1990hierarchy,sherali1994hierarchy,sherali2009reformulation}.

RLT has been shown to provide strong relaxations~\cite{sherali1992new,sherali2000reduced}, but
this comes at the cost of excessive numbers of cuts.
To address this, Sherali and Tuncbilek~\cite{sherali1997new} proposed a technique to
add a subset of RLT cuts, depending on signs of coefficients of monomial terms in the
original constraints and the RLT constraints.
Furthermore, the reduced 
RLT technique for polynomial problems
containing linear equality constraints~\cite{liberti2004reduction,liberti2004reformulation,liberti2005linearity,sherali2012reduced}
enforces a subset of polynomial relations via RLT cuts and discards the original polynomial constraints, and it is shown to yield
equivalent representations with fewer nonlinear terms.
For binary quadratic problems containing linear equality constraints, a compact linearization is
introduced through a selective application of RLT~\cite{liberti2007compact,mallach2018compact},
which was generalized to binary quadratic problems with arbitrary linear constraints in the inductive
linearization approach developed by Mallach~\cite{mallach2021inductive}.

We focus on RLT for bilinear \firstrevision{terms}, which is of particular interest due
to the numerous applications whose models involve nonconvex quadratic
nonlinearities~\cite{adams1993mixed,frank2012optimal,frank2012optimal2,misener2009advances,buchheim2010exact,castillo2005optimization}.
Even in the bilinear case, large numbers of factors to be
multiplied and of RLT cuts that are generated as a result remain an issue that can lead to considerable slowdowns, both due
to the cost of cut separation and the large sizes of resulting LP relaxations.

The first contribution of this paper is a new approach to applying RLT to
MILPs.
Unlike the approaches that only introduce multilinear relations via
multiplication~\cite{sherali1990hierarchy,sherali1994hierarchy},
this approach detects and enforces bilinear relations that are already implicitly present in the
model.
A bilinear product relation where one multiplier is a binary variable and the other
multiplier is a variable with finite bounds can be equivalently written as two linear constraints.
We identify such pairs of linear constraints
that implicitly encode a bilinear product relation,
then utilize this relation in the generation of RLT cuts.
We present two strategies for filtering implicit product relations and compare these strategies
experimentally.

The second contribution of this paper addresses the major bottleneck for applying RLT successfully
in practice, which stems from the prohibitive cost of separating all RLT cuts,
by proposing an efficient separation algorithm.
This algorithm considers the signs of bilinear relation violations in a current LP
relaxation solution and the signs of coefficients in linear constraints in order to ignore
combinations of factors that will not produce a violated inequality.
Furthermore, we propose a technique which projects the linear constraints onto a reduced space
and constructs RLT cuts based on the resulting much smaller system.
Depending on the implementation, this projection may have an added effect of filtering cuts.

The rest of the paper is organized as follows.
In Section~\ref{section:background}, RLT for bilinear \firstrevision{terms} is explained.
In Section~\ref{section:detection}, we describe the technique for deriving bilinear product relations from
MILP constraints and discuss product relation filtering strategies.
Section~\ref{section:separation} presents the new cut separation algorithm, and
computational results are presented in Section~\ref{section:results}.
Section~\ref{section:conclusion} concludes the paper.

\section{RLT for Bilinear \firstrevision{Terms}}\label{section:background}

We consider mixed-\firstrevision{binary} nonlinear programs of the extended form where auxiliary
variables $w$ are introduced for all bilinear \firstrevision{terms}:

\begin{subequations}
\label{minlp-problem}
\begin{align}
\min\; &\mathbf{c}^{\text{T}}\mathbf{x}\\
\text{s.t.}\; &A\mathbf{x} \leq \mathbf{b},\label{eq:linear}\\
&g(\mathbf{x},\mathbf{w}) \leq 0,\\
&x_ix_j \lesseqgtr w_{ij} \text{ for all } (i,j) \in \mathcal{I}^w, \label{productrels}\\ 
&\mathbf{\underline{x}} \leq \mathbf{x} \leq \mathbf{\overline{x}}, ~\mathbf{\underline{w}} \leq \mathbf{w} \leq \mathbf{\overline{w}},\\
&x_j \in \mathbb{R} \text{ for all } j \in \mathcal{I}^c, ~x_j \in \{0,1\} \text{ for all } j \in \mathcal{I}^b,
\end{align}
\end{subequations}
with $\mathcal{I} = \mathcal{I}^c \cup \mathcal{I}^b$ being a disjoint partition of variables $\mathbf{x}$
and $\mathbf{x}$ having dimension $|\mathcal{I}| = n$.
In the above formulation,
$\mathbf{\underline{x}}$, $\mathbf{\overline{x}} \in \overline{\mathbb{R}}^n$, $\mathbf{\underline{w}}, \mathbf{\overline{w}} \in \overline{\mathbb{R}}^{|\mathcal{I}^w|}$
($\overline{\mathbb{R}} = \mathbb{R} \cup \{-\infty,+\infty\}$), $\mathbf{c} \in \mathbb{R}^n$
and $\mathbf{b} \in \mathbb{R}^{m^{(l)}}$ are constant vectors and
$A \in \mathbb{R}^{m^{(l)} \times n}$
is a coefficient matrix, and
the function $g$
defines the nonlinear constraints.
Constraint~\eqref{productrels} defines the bilinear product relations in the problem and allows for
inequalities and equations.
Let $\mathcal{I}^p$ denote the set of indices of all variables that participate in bilinear
product relations~\eqref{productrels}, that is, all indices $i$ such that there exists an index $j$ for which
$(i,j) \in \mathcal{I}^w$.

\firstrevision{The techniques presented in this paper are also applicable to problems with general integer variables, i.e., to general mixed-integer nonlinear programs (MINLPs),
but general integer variables are treated similarly to continuous variables in our methods.
Therefore, for simplicity of presentation, we do not include general integer variables in the above
problem formulation.}

Solvers typically employ McCormick inequalities~\cite{mccormick1976computability} to
construct an LP relaxation of constraints~\eqref{productrels}.
These inequalities describe
the convex hull of the set given by the relation $x_ix_j \lesseqgtr w_{ij}$:
\begin{subequations}\label{eq:mccormickgeq}
\begin{align}
  \overline{x}_ix_j + x_i\underline{x}_j - \overline{x}_i\underline{x}_j &\geq w_{ij},\\
  \underline{x}_ix_j + x_i\overline{x}_j - \underline{x}_i\overline{x}_j &\geq w_{ij},
\end{align}
\end{subequations}
as relaxation of $x_ix_j \geq w_{ij}$, and
\begin{subequations}\label{eq:mccormickleq}
\begin{align}
  \overline{x}_ix_j + x_i\overline{x}_j - \overline{x}_i\overline{x}_j &\leq w_{ij},\\
  \underline{x}_ix_j + x_i\underline{x}_j - \underline{x}_i\underline{x}_j &\leq w_{ij},
\end{align}
\end{subequations}
as relaxation of $x_ix_j \leq w_{ij}$.

In the presence of linear constraints~\eqref{eq:linear}, this relaxation can be strengthened by
adding RLT cuts.
RLT is based on multiplying constraints by nonnegative factors in order to obtain valid
nonlinear inequalities, and then linearizing the result by substituting the nonlinear terms, where
possible, with equivalent linear terms, and outer approximating the remaining nonlinear terms.

Consider a linear constraint:
$\sum_{k=1}^n a_{1k}x_k \leq b_1.$
Multiplying this constraint by nonnegative bound factors $(x_j-\underline{x}_j)$ and $(\overline{x}_j-x_j)$,
where $\underline{x}_j$ and $\overline{x}_j$ are finite, yields
valid nonlinear inequalities.
\firstrevision{In the RLT version which we consider in this paper, a constraint is multiplied by
one bound factor at a time.}
We will derive the RLT cut using the lower bound factor.
The derivation is analogous for the upper bound factor.
The multiplication, referred to as the reformulation step, yields:
$$\sum_{k=1}^n a_{1k}x_k(x_j-\underline{x}_j) \leq b_1(x_j-\underline{x}_j).$$
%
This nonlinear inequality is then linearized in order to obtain a valid linear inequality.
%
%
%
%
The following linearizations are applied to each nonlinear term $x_kx_j, k=1,\dots,n$:
\begin{itemize}
\item $x_kx_j$ is replaced by $w_{kj}$ if the relation $x_kx_j \leq w_{kj}$ exists in the problem
      and $a_{1k} \leq 0$, or if the relation $x_kx_j \geq w_{kj}$ exists and $a_{1k} \geq 0$,
      or if the relation $x_kx_j = w_{kj}$ exists in the problem,
\item if $k=j \in \mathcal{I}^{b}$, then $x_kx_j = x_j^2 = x_j$,
\item if $k=j \notin \mathcal{I}^b$, then $x_kx_j = x_j^2$ is outer approximated by a secant from above or by a tangent
from below, depending on the sign of the coefficient,
\item if $k \neq j$, $k,j \in \mathcal{I}^{b}$ and one of the four clique constraints is implied
  by the linear constraints~\eqref{eq:linear}, then:
  $x_k + x_j \leq 1 \Rightarrow x_kx_j = 0$; $x_k - x_j \leq 0 \Rightarrow x_kx_j = x_k$;
  $-x_k + x_j \leq 0 \Rightarrow x_kx_j = x_j$; $-x_k - x_j \leq -1 \Rightarrow x_kx_j = x_j + x_j - 1$,
\item otherwise, $x_kx_j$ is replaced by its McCormick relaxation.

\end{itemize}

In this way we obtain a valid linear inequality.
The key step is the replacement of products $x_kx_j$ with the variables $w_{kj}$.
When a bilinear product relation $x_kx_j \lesseqgtr w_{kj}$ does not hold for the current relaxation
solution, this substitution may lead to an increase in the violation of the inequality, thus
possibly producing a cut that is violated by the relaxation solution.

In the case that we have a linear equation constraint $\sum_{k=1}^n a_{1k}x_k = b_1$ \emph{and} all nonlinear terms
can be replaced using equality relations, then RLT produces an equation cut.
Otherwise, the equation constraint is treated as two inequalities $\sum_{k=1}^n a_{1k}x_k \leq b_1$ and
$\sum_{k=1}^n a_{1k}x_k \geq b_1$ to produce inequality cuts.
Note that nonnegative constraint factors $(b_i - \mathbf{a}_i\mathbf{x})$ could also be used as multipliers
similarly to bound factors, but we do not do so in our implementation. 

\section{Detection of Implicit Products}\label{section:detection}

This section explains how bilinear product relations can be derived from mixed-integer linear
constraints.

\subsection{Deriving Product Relations from MILP Constraints}\label{subsection:detection_theory}

Let us begin with the reverse direction:
Consider a product relation
$w_{ij} = x_ix_j$, where $x_i$ is binary.
It can be equivalently rewritten as two implications: $x_i = 1 \Rightarrow w_{ij} = x_j$ and
$x_i = 0 \Rightarrow w_{ij} = 0$ .
With the use of the big-M technique, these implications can be represented as linear constraints,
provided that the bounds of $x_j$ are finite:
\begin{subequations}\label{eq:bigm_leq}
\begin{align}
\hspace*{1.7em} w_{ij} - x_j - \underline{x}_jx_i &\leq -\underline{x}_j,\\
w_{ij} - \overline{x}_jx_i &\leq 0
\end{align}
\end{subequations}
to ensure that $w_{ij} \leq x_ix_j$ holds, and
\begin{subequations}\label{eq:bigm_geq}
\begin{align}
-w_{ij} + x_j + \overline{x}_jx_i &\leq \overline{x}_j,\\
-w_{ij} + \underline{x}_jx_i &\leq 0
\end{align}
\end{subequations}
to ensure that $w_{ij} \geq x_ix_j$ holds.

Conversely, linear constraints with binary variables can be analyzed in order to
detect constraint pairs of the forms~\eqref{eq:bigm_leq} $\Leftrightarrow w_{ij} \leq x_ix_j$ or~\eqref{eq:bigm_geq} $\Leftrightarrow w_{ij} \geq x_ix_j$.
This method can be generalized to derive bilinear relations of the form
\begin{equation}\label{eq:implicit_product}Ax_i + Bw_{ij} + Cx_j + D \lesseqgtr x_ix_j\end{equation}
with $A,B,C,D \in \mathbb{R}$.
The process of deriving a bilinear product relation from two suitable linear constraints is
described in the following theorem.

\begin{theorem}\label{thm:hiddenproduct}
Consider two linear constraints depending on the same three variables $x_i$, $x_j$ and $w_{ij}$,
where $x_i$ is binary:
\begin{subequations}\label{eq:detection_linrels}
\begin{gather}
\label{eq:detection_linrels1}
a_1x_i + b_1w_{ij} + c_1x_j \leq d_1,\\
\label{eq:detection_linrels2}
a_2x_i + b_2w_{ij} + c_2x_j \leq d_2.
\end{gather}
\end{subequations}
If $b_1b_2 > 0$ and $\gamma = c_2b_1 - b_2c_1 \neq 0$, then these constraints imply the
product relation
\begin{align}
\label{eq:implicit_product_formula}
\frac{b_2(a_1 - d_1) + b_1d_2}{\gamma}x_i + \frac{b_1b_2}{\gamma}w_{ij} + \frac{b_1c_2}{\gamma}x_j - \frac{b_1d_2}{\gamma}
\begin{cases}
  \leq x_ix_j & \text{ if $b_1\gamma > 0$,} \\
  \geq x_ix_j & \text{ if $b_1\gamma < 0$.}
\end{cases}
\end{align}
\end{theorem}
\begin{proof}
We begin by writing the
bilinear relation \eqref{eq:implicit_product}, treating its coefficients and inequality sign as
unknown,
and reformulating it as two implications:
\begin{subequations}\label{eq:detection_implics1}
\begin{align}
&x_i = 1 &~\Rightarrow~ &Bw_{ij} + (C-1)x_j \lesseqgtr -D - A,\\
&x_i = 0 &~\Rightarrow~ &Bw_{ij} + Cx_j \lesseqgtr -D,
\end{align}
\end{subequations}
where the inequality sign must be identical in both implied inequalities.
Similarly, we rewrite constraints~\eqref{eq:detection_linrels} with scaling parameters $\alpha$ and $\beta$
as implications:
\begin{subequations}\label{eq:detection_implics2}
\begin{align}
&x_i = 1 &~\Rightarrow~ &\alpha b_1w_{ij} + \alpha c_1x_j \lesseqgtr \alpha(d_1 - a_1),\\
&x_i = 0 &~\Rightarrow~ &\beta b_2w_{ij} + \beta c_2x_j \lesseqgtr \beta d_2,
\end{align}
\end{subequations}
where the inequality signs depend on the signs of $\alpha$ and $\beta$.

The goal is to find the coefficients $A,B,C$ and $D$ and the inequality sign of the
relation~\eqref{eq:implicit_product} that can be derived from~\eqref{eq:detection_implics2}
for suitable $\alpha$ and $\beta$.
We require that coefficients and inequality signs in implications~\eqref{eq:detection_implics1}
and~\eqref{eq:detection_implics2} match:
%
\begin{align*}
B &= \alpha b_1
&
B &= \beta b_2\\
C - 1 &= \alpha c_1
&
C &= \beta c_2
&
\alpha\beta > 0\\
D + A &= \alpha (a_1-d_1)
&
D &= -\beta d_2
\end{align*}
%
The equations $\alpha b_1-\beta b_2=0$ and $\alpha c_1-\beta c_2=-1$ determine the scaling parameters $\alpha$ and $\beta$ as
$$\alpha = \frac{b_2}{c_2b_1 - b_2c_1}, ~\beta = \frac{b_1}{c_2b_1 - b_2c_1}.$$
Therefore we must require that $c_2b_1 - b_2c_1 \neq 0$.
Since the inequality signs must be similar in the two scaled linear relations, the signs of $\alpha$
and $\beta$ must be equal.
If $\alpha$ and $\beta$ are positive, all inequalities are with ``less than or equal''.
Otherwise, they become ``greater than or equal''.

Solving the resulting system yields:
\begin{gather*}
b_1b_2 > 0, ~A = (1/\gamma)(b_2(a_1 - d_1) + b_1d_2)\\
B = b_1b_2/\gamma, ~C = b_1c_2/\gamma, ~D = -b_1d_2/\gamma, ~\gamma \neq 0,
\end{gather*}
where $\gamma = c_2b_1 - b_2c_1$
and the inequality sign is `$\leq$` if $b_1\gamma \geq 0$, and `$\geq$` if $b_1\gamma \leq 0$.
Thus, the bilinear relation is obtained:
$$Ax_i + Bw_{ij} + Cx_j + D \leq x_ix_j ~\text{ if } \frac{b_1}{c_2b_1 - b_2c_1} \geq 0,$$
$$Ax_i + Bw_{ij} + Cx_j + D \geq x_ix_j ~\text{ if } \frac{b_1}{c_2b_1 - b_2c_1} \leq 0.$$
\end{proof}

Although the conditions of the theorem are sufficient for the bilinear product relation to be
implied by the linear constraints, in practice we may check more conditions before deriving such a
relation.
In particular:

\begin{enumerate}
\item At least one of the coefficients $a_1$ and $a_2$ must be nonzero.
Otherwise, the product relation is always implied by the linear constraints,
including when $0 < x_i < 1$.
In our implementation, this condition is always enforced.
\item The signs of the coefficients of the binary variable $x_i$ must be different, that is, one
linear relation is more restrictive when $x_i = 1$ and the other when $x_i = 0$.
While this is not necessary for the non-redundancy of the derived product relation, by requiring
this we focus on stronger implications (for instance, for a linear relation $a_1x_i + b_1w_{ij} +
c_1x_j \leq d_1$ with $a_1 > 0$, we use the more restrictive implication $x_i = 1 ~\Rightarrow
b_1w_{ij} + c_1x_j \leq d_1 - a_1$ rather than the less restrictive implication $x_i = 0
~\Rightarrow b_1w_{ij} + c_1x_j \leq d_1$).
In our implementation, a parameter controls whether this condition is enforced.
This condition is enabled by default.
\end{enumerate}

In separation, the product relation \eqref{eq:implicit_product} is treated similarly to product
relations $w_{ij} \lesseqgtr x_ix_j$, with the linear left-hand side $Ax_i + Bw_{ij} + Cx_j + D$
being used instead of the individual auxiliary variable $w_{ij}$.

\subsection{Detection Algorithm}\label{detection_algorithm}

The detection algorithm searches for suitable pairs of linear relations and derives product
relations from them.
Let $x_i$, as before, be a binary variable.
The following relation types are considered as candidates for the first relation in such a pair:
\begin{itemize}
\item implied relations of the form $x_i = \xi ~\Rightarrow~ \tilde b_1w_{ij} + \tilde c_1x_j \leq
\tilde d_1$, where $\xi = 0$ or $\xi = 1$; such relations are derived from linear constraints and
stored in a hashtable which maps
variable triples to arrays of linear constraints that contain only the variables in the triple;
\item implied bounds of the form $x_i = \xi ~\Rightarrow~ w_{ij} \leq \tilde d_1$, stored as sorted
arrays for each $x_i$.
\end{itemize}

The second relation in a pair can be:
\begin{itemize}
\item an implied relation of the form $x_i = \overline{\xi} ~\Rightarrow~ \tilde b_2w_{ij} + \tilde
c_2x_j \leq \tilde d_2$, where $\overline{\xi}$ is the complement of $\xi$;
\item if $w_{ij}$ is non-binary, an implied bound of the form $x_i = \overline{\xi} ~\Rightarrow~
w_{ij} \leq \tilde d_2$;
\item if $w_{ij}$ is binary, a clique containing the complement of $x_i$ if $\xi = 1$ or $x_i$ if
$\xi = 0$, and $w_{ij}$ or its complement; Cliques are constraints of the form: $\sum_{k \in
\mathcal{J}}x_k + \sum_{k \in \overline{\mathcal{J}}}(1-x_k) \leq 1$, where $\mathcal{J} \subseteq
\mathcal{I}^b$, $\overline{\mathcal{J}} \subseteq \mathcal{I}^b$ and $\mathcal{J} \cap \overline{\mathcal{J}} = \emptyset$. Such constraints are stored
in a clique table;
\item a constraint on $x_j$ and $w_{ij}$, which can be stored as
an implied bound, clique or general linear constraint; the latter are stored
in a hastable similar to the three-variable implied relations;
\item  a global bound on $w_{ij}$.
%
\end{itemize}
If the first relation is an implied bound $x_i = \xi ~\Rightarrow~ w_{ij} \leq \tilde d_1$, then
only relations between $w_{ij}$ and one other, arbitrarily chosen, variable are considered.

Once a suitable pair of constraints is found, the bilinear product relation is derived from it as
described in Theorem~\ref{thm:hiddenproduct}.
All bilinear \firstrevision{terms}, whether implicitly or explicitly defined, are stored with the use of a
specialized structure that contains the information on the product type, the so-called uplocks and
downlocks~\cite{Achterberg2007a}, which indicate whether increasing or decreasing the value of the product expression can
make some constraints in the problem infeasible, and arrays of variables and linear expressions
that the product is linked to via a bilinear product relation.

\subsection{Identifying redundant implicit relations}\label{sec:redchecks}

To determine when a product relation~\eqref{eq:implicit_product_formula} is redundant with
respect to at least one of the linear relations~\eqref{eq:detection_linrels}, we first recall that
$b_1b_2 > 0$.
Consider the case $b_1 > 0$ and $b_2 > 0$ and rewrite the inequalities as:
\begin{gather}
\label{eq:detection_linrels1p}\tag{\ref{eq:detection_linrels1}'}
w_{ij} \leq d_1/b_1 - (a_1/b_1)x_i - (c_1/b_1)x_j,\\
\label{eq:detection_linrels2p}\tag{\ref{eq:detection_linrels2}'}
w_{ij} \leq d_2/b_2 - (a_2/b_2)x_i - (c_2/b_2)x_j,\\
\label{eq:implicit_product_formulap}\tag{\ref{eq:implicit_product_formula}'}
w_{ij} \leq \frac{\gamma}{b_1b_2}x_ix_j - \left(\frac{a_1 - d_1}{b_1} + \frac{d_2}{b_2}\right)x_i - \frac{c_2}{b_2}x_j + \frac{d_2}{b_2}.
\end{gather}
The product relation \eqref{eq:implicit_product_formulap} is redundant with respect to \eqref{eq:detection_linrels1p} if: 
$$b_2d_1 - a_1b_2x_i - b_2c_1x_j \leq \gamma x_ix_j - (b_2(a_1 - d_1) + b_1d_2)x_i - b_1c_2x_j + b_1d_2 \quad \Leftrightarrow$$
$$b_2d_1 - b_2c_1x_j \leq \gamma x_ix_j - (b_1d_2 - b_2d_1)x_i - b_1c_2x_j + b_1d_2 \quad \Leftrightarrow$$
$$(b_2d_1 - b_1d_2)(1 - x_i) + \gamma x_j \leq \gamma x_ix_j  \quad \Leftrightarrow$$
$$(b_2d_1 - b_1d_2)(1 - x_i) + \gamma x_j(1 - x_i) \leq 0 \quad \Leftrightarrow$$
$$\gamma x_j \leq b_1d_2 - b_2d_1.$$
It is redundant with respect to \eqref{eq:detection_linrels2p} if: 
$$b_1d_2 - a_2b_1x_i - b_1c_2x_j \leq \gamma x_ix_j - (b_2(a_1 - d_1) + b_1d_2)x_i - b_1c_2x_j + b_1d_2 \quad \Leftrightarrow$$
$$(b_2(a_1 - d_1) + b_1(d_2 - a_2))x_i \leq \gamma x_ix_j \quad \Leftrightarrow$$
$$b_1d_2 - b_2d_1 + a_1b_2 - a_2b_1 \leq \gamma x_j.$$
For the case of $b_1<0$ and $b_2<0$, we have:
\begin{gather*}
\label{eq:detection_linrels1pp}\tag{\ref{eq:detection_linrels1}''}
w_{ij} \geq d_1/b_1 - (a_1/b_1)x_i - (c_1/b_1)x_j,\\
\label{eq:detection_linrels2pp}\tag{\ref{eq:detection_linrels2}''}
w_{ij} \geq d_2/b_2 - (a_2/b_2)x_i - (c_2/b_2)x_j,\\
\label{eq:implicit_product_formulapp}\tag{\ref{eq:implicit_product_formula}''}
w_{ij} \geq \frac{\gamma}{b_1b_2}x_ix_j - \left(\frac{a_1 - d_1}{b_1} + \frac{d_2}{b_2}\right)x_i - \frac{c_2}{b_2}x_j + \frac{d_2}{b_2},
\end{gather*}
and the product relation \eqref{eq:implicit_product_formulapp} is redundant with respect to \eqref{eq:detection_linrels1pp} if: 
$$\gamma x_j \geq b_1d_2 - b_2d_1,$$
and with respect to \eqref{eq:detection_linrels2pp} if: 
$$b_1d_2 - b_2d_1 + a_1b_2 - a_2b_1 \geq \gamma x_j.$$
All in all, the bilinear product relation is nonredundant if:
$$
\begin{cases}
x_j \in (1/\gamma(b_1d_2 - b_2d_1 + a_1b_2 - a_2b_1), 1/\gamma(b_1d_2 - b_2d_1)) \text{ if } b_1\gamma > 0,\\
x_j \in (1/\gamma(b_1d_2 - b_2d_1), 1/\gamma(b_1d_2 - b_2d_1 + a_1b_2 - a_2b_1)) \text{ if } b_1\gamma < 0.
\end{cases}
$$
The detection algorithm only saves products for which the intersection of this interval
with the domain of $x_j$ is nonempty.

\section{Separation Algorithm}\label{section:separation}

We present a new algorithm for separating RLT cuts within an LP-based branch-and-bound
solver.
The branch-and-bound algorithm builds LP relaxations of problem~\eqref{minlp-problem}
by constructing linear underestimators of functions $g$ in the constraint $g(\mathbf x,\mathbf w) \leq 0$
and McCormick inequalities for constraints~\eqref{productrels}.

Let $(\mathbf{x}^*,\mathbf{w}^*)$ be the solution of an LP relaxation at a node of the branch-and-bound tree, and
suppose that $(\mathbf{x}^*,\mathbf{w}^*)$ violates the relation $x_ix_j \lesseqgtr w_{ij}$ for some
$i,j \in \mathcal{I}^w$.
Separation algorithms generate cuts that separate $(\mathbf{x}^*,\mathbf{w}^*)$ from the feasible region, and
add those cuts to the solver's cut storage.

There are two types of cut storages: the cuts in the storage of the current node are selected by
the solver's internal algorithms and added or discarded immediately, and cuts in the cut pool are
retained even if the current LP solution does not violate them, and such cuts may be added in other
nodes of the branch-and-bound tree.

\subsection{Standard Separation}\label{sec:simplesepa}

The standard separation algorithm, which will serve as a baseline for
comparisons,
iterates over all linear constraints.
For each constraint, it iterates over all variables $x_j$ that participate in bilinear relations
and for each such variable, the constraint is multiplied by bound factors $(x_j - \underline{x}_j)$ and
$(\overline{x}_j - x_j)$.
RLT cuts are then obtained by linearizing the result of the multiplication, and 
violated cuts are
added to the MINLP solver's cut storage.

\subsection{Row Marking}\label{subsection:marking}

Let the bound factors be denoted as $f_j^{(\ell)}(\mathbf{x}) = x_j-\underline{x}_j$ and $f_j^{(u)}(\mathbf{x}) = \overline{x}_j - x_j$.
Consider a linear constraint multiplied by a bound factor:
\begin{equation}\label{eq:sepa}f_j^{(.)}(\mathbf{x})\mathbf{a}_r\mathbf{x} \leq f_j^{(.)}(\mathbf{x})b_r.
\end{equation}

The $i$th nonlinear term is $a'_{ri}x_ix_j$, where $a'_{ri} = a_{ri}$ when
multiplying by $(x_j - \underline{x}_j)$ and $a'_{ri} = -a_{ri}$ when multiplying by $(\overline{x}_j - x_j)$.
Following the procedure described in Section~\ref{section:background}, RLT may replace the product $x_ix_j$ with $w_{ij}$.
The product can also be replaced with a linear expression, but this does not change
the reasoning, and we will only use $w_{ij}$ in this section.

If $w^*_{ij} \neq x^*_ix^*_j$, then such a
replacement will change the violation of~\eqref{eq:sepa}.
The terms whose replacement will increase the violation are of interest, that is, the terms
where:
$$a'_{ri}x^*_ix^*_j \leq a'_{ri}w^*_{ij}.$$

This determines the choice of bound factors to multiply with:
$$x^*_ix^*_j < w^*_{ij} ~\Rightarrow~ \begin{array}{c}\text{multiply by }(x_j - \underline{x}_j) \text{ if }a_{ri} > 0,\\
                                                      \text{multiply by }(\overline{x}_j - x_j) \text{ if }a_{ri} < 0,\end{array}$$
$$x^*_ix^*_j > w^*_{ij} ~\Rightarrow~ \begin{array}{c}\text{multiply by }(\overline{x}_j - x_j) \text{ if }a_{ri} > 0,\\
                                                      \text{multiply by }(x_j - \underline{x}_j) \text{ if }a_{ri} < 0.\end{array}$$

The separation algorithm is initialized by creating the following data structures:
\begin{itemize}
\item sparse arrays which allow to access all variables appearing in bilinear \firstrevision{terms} together
with a given variable,
\item a hashtable for accessing the bilinear \firstrevision{term} by two given variables.
\end{itemize}

For each variable $x_i$, linear rows are marked in order to inform the separation algorithms which bound
factors of $x_i$ they should be multiplied with, if any.
The algorithm can work with inequality rows in both `$\leq$' and `$\geq$' forms as
well as equation rows.
For each bilinear \firstrevision{term} $x_ix_j$, the row marking algorithm iterates over all linear rows that
contain $x_j$ with a nonzero coefficient.
These rows are stored in a sparse array and have one of the following marks:
\begin{itemize}
\item MARK\_LT: the row contains a term $a_{rj}x_j$ such that $a_{rj}x^*_ix^*_j < a_{rj}w^*_{ij}$;
\item MARK\_GT: the row contains a term $a_{rj}x_j$ such that $a_{rj}x^*_ix^*_j > a_{rj}w^*_{ij}$;
\item MARK\_BOTH: the row contains terms fitting both cases above.
\end{itemize}
Row marks are stored in two sparse arrays, \textit{row\_idcs} and
\textit{row\_marks}, the first storing sorted row indices and the second storing the corresponding
marks.
In the algorithm below, we use the notation \firstrevision{\textit{marks(r)}} to denote accessing the mark of row
$r$ by performing a search in \textit{row\_idcs} and retrieving the corresponding entry in
\textit{row\_marks}, \firstrevision{or adding such an entry if it does not exist}.
\firstrevision{In the implementation of the algorithm, marks are represented by bit flags.}
We also define a sparse matrix $W$ with entries $w_{ij}$.

\firstrevision{Algorithm~\ref{alg:marking} shows} the pseudocode for the row marking algorithm.
%
%

The \firstrevision{marking-based separation} algorithm iterates over the sparse array of marked rows and generates RLT cuts for
the following combinations of linear rows and bound factors:

\begin{itemize}
\item If \firstrevision{the mark is set to} MARK\_LT, then ``$\leq$'' constraints are multiplied with the lower bound
factor and ``$\geq$'' constraints are multiplied with the upper bound factor;
\item If \firstrevision{the mark is set to} MARK\_GT, then ``$\leq$'' constraints are multiplied with the upper bound
factor and ``$\geq$'' constraints are multiplied with the lower bound factor;
\item If \firstrevision{the mark is set to} MARK\_BOTH, then both ``$\leq$'' and ``$\geq$'' constraints are multiplied
with both the lower and the upper bound factors;
\item Marked equality constraints are always multiplied with $x_i$ itself.
\end{itemize}

\begin{algorithm}\caption{Row marking}\label{alg:marking}
\begin{algorithmic}
\State \textbf{Input: } $x^*, w^*, W$
\State $marks \coloneqq \emptyset$
\For{$i \in \mathcal{I}^p, j \in nnz(\bm w_i)$}
  \For{$r$ such that $j \in nnz(\bm a_r)$}
    \If{$r \notin \firstrevision{row\_idcs}$}
      \State \firstrevision{set} $marks(r)$ \firstrevision{to} $0$
    \EndIf
    \If{$a_{rj}x^*_ix^*_j < a_{rj}w^*_{ij}$}
      \State \firstrevision{set} $marks(r)$ \firstrevision{to}
      \firstrevision{$\begin{cases}\text{MARK\_LT if }marks(r)=0,\\ \text{MARK\_BOTH if }marks(r)=\text{MARK\_GT}\end{cases}$}
    \ElsIf{\firstrevision{$a_{rj}x^*_ix^*_j > a_{rj}w^*_{ij}$}}
      \State \firstrevision{set} $marks(r)$ \firstrevision{to}
      \firstrevision{$\begin{cases}\text{MARK\_GT if }marks(r)=0,\\ \text{MARK\_BOTH if }marks(r)=\text{MARK\_LT}\end{cases}$}
    \EndIf
  \EndFor
\EndFor
\end{algorithmic}
\end{algorithm}

\subsection{Projection Filtering}\label{subsection:projection}

\firstrevision{Similarly to row marking, the projection filtering technique aims to avoid
the processing of factors that will not yield a violated cut.
The key observation is that if a product variable is at its bound, then McCormick inequalities
for a product relation involving this variable imply the product relation.
Projection filtering projects the linear constraints onto the space of variables for which
$\underline{x}_j < x^*_j < \overline{x}_j$.
Violations of RLT cuts generated from projected inequalities are checked
and the factor combinations for which such cuts are not violated are discarded in the current separation round.}

\firstrevision{If all McCormick inequalities are present in the LP relaxation, then
the RLT cut generated from an inequality is violated if and only if the RLT cut
generated from its projected version is violated.
Since $\mathbf x^*$ is a basic solution, in practice many variables tend to have
values at a bound, and the projected system often has a considerably smaller size than the original system.
Thus, projection filtering can reduce the computational cost of separation.}

\firstrevision{Depending on the solver, however, the McCormick inequalities may not all be satisfied at an LP solution.
This is the case in the MINLP solver SCIP~\cite{BestuzhevaBesanconEtal2023_Enabling,bestuzheva2023global}, which dynamically separates McCormick cuts.
The violations of RLT cuts generated from original and projected inequalities are then not guaranteed to be
equal.
In this case, projection filtering has an additional effect: for violated bilinear product relations involving
variables whose values are at bound, the violation of the product will be disregarded when checking
the violation of RLT cuts.
Thus, adding McCormick cuts will be prioritized over adding RLT cuts, potentially reducing the
number of RLT cuts added to the LP relaxation.}

\firstrevision{In the remainder of this section, we discuss these properties and their implications formally.}
Recall the McCormick inequalities~\eqref{eq:mccormickgeq} and~\eqref{eq:mccormickleq}:
\begin{subequations}
\begin{align}
&\overline{x}_ix_j + x_i\underline{x}_j - \overline{x}_i\underline{x}_j \geq w_{ij},\label{mc1}\\
&\underline{x}_ix_j + x_i\overline{x}_j - \underline{x}_i\overline{x}_j \geq w_{ij},\label{mc2}\\
&\overline{x}_ix_j + x_i\overline{x}_j - \overline{x}_i\overline{x}_j \leq w_{ij},\label{mc3}\\
&\underline{x}_ix_j + x_i\underline{x}_j - \underline{x}_i\underline{x}_j \leq w_{ij}.\label{mc4}
\end{align}
\end{subequations}
If at least one of the variables $x_i$ and $x_j$ has a value equal to one of its
bounds, then the McCormick relaxation is tight for the relation $w_{ij} = x_ix_j$.
\firstrevision{This can be seen directly by substituting the variable with its bound value, which leads to
two of the above inequalities together being equivalent to the product relation:}
\firstrevision{
\begin{itemize}
\item if $x_i = \underline{x}_i$, then \{\eqref{mc2}, \eqref{mc4}\} $\Leftrightarrow~ x_ix_j = w_{ij}$;
\item if $x_i = \overline{x}_i$, then \{\eqref{mc1}, \eqref{mc3}\} $\Leftrightarrow~ x_ix_j = w_{ij}$;
\item if $x_j = \underline{x}_j$, then \{\eqref{mc1}, \eqref{mc4}\} $\Leftrightarrow~ x_ix_j = w_{ij}$;
\item if $x_i = \overline{x}_i$, then \{\eqref{mc2}, \eqref{mc3}\} $\Leftrightarrow~ x_ix_j = w_{ij}$.
\end{itemize}
}
Therefore, if $x_i$ or $x_j$ is at a bound and the McCormick inequalities are satisfied, then the
product relation is also satisfied.
We will describe the equality case here, and the reasoning is analogous for the inequality cases of $x_ix_j \lesseqgtr w_{ij}$.

Consider the linear system $A\mathbf x \leq \mathbf b$ projected onto the set of variables whose values
are not equal to either of their bounds:
$$\sum_{k \in \mathcal{J}^1} a_{rk}x_k \leq b_r - \sum_{k \in \mathcal{J}^2} a_{rk}x^*_k, ~\forall r \in 1,\dots,m^{(l)},$$
where $\mathcal{J}^1 \subseteq \mathcal{I}$ is the set of all problem variables whose values in the
solution $\mathbf{x}^*$ of the current LP relaxation are not equal to one of their bounds, and
$\mathcal{J}^2 = \mathcal{I} \setminus \mathcal{J}^1$.

Violation is then first checked for RLT cuts generated based on the projected linear system.
Only if such a cut, which we will refer to as a projected RLT cut, is violated, then the RLT cut
for the same bound factor and the corresponding constraint in the original linear system will be
constructed.

\firstrevision{}In the projected system multiplied with a bound factor $f_j^{(.)}(\mathbf{x})$:
$$f_j^{(.)}(\mathbf{x})\cdot\sum_{k \in \mathcal{J}^1} a_{rk}x_k \leq f_j^{(.)}(\mathbf{x})(b_r - \sum_{k \in \mathcal{J}^2} a_{rk}x^*_k), ~\forall r \in 1,\dots,m^{(l)},$$
the only nonlinear terms are $x_jx_k$ with $k \in \mathcal{J}^1$,
and therefore, no
substitution $x_ix_k \rightarrow w_{ik}$ is performed for $k \in \mathcal{J}^2$.
\firstrevision{The implications depend on whether McCormick inequalities hold:}
\begin{itemize}
\item If the McCormick inequalities for $x_i$, $x_k$, and $w_{ik}$ hold, then $x^*_ix^*_k = w^*_{ik}$
for $k \in \mathcal{J}^2$, and \firstrevision{if substitution was performed, it} would not change the violation of the inequality.
Therefore, \firstrevision{in this case} checking the violation of a projected RLT cut
is equivalent to checking the violation of a full RLT cut.
\item \firstrevision{If the McCormick inequalities for $x_i$, $x_k$, and $w_{ik}$ do not hold, then it is possible that $x^*_ix^*_k \neq w^*_{ik}$ for some $k \in \mathcal{J}^2$, but these violations
will not contribute to the violation of the projected RLT cut, since no
substitution is performed for fixed variables.
Thus, the violations of products involving variables at bound may be ignored, leading to the aforementioned prioritization of McCormick inequalities.}
\end{itemize}

%
%
%
%

\section{Computational Results}\label{section:results}

\subsection{Setup}

We tested the proposed methods on the
MINLPLib\footnote{\url{https://www.minlplib.org}}~\cite{10.1287/ijoc.15.1.114.15159} test set and a
test set comprised of instances from MIPLIB3, MIPLIB 2003, 2010 and 2017~\cite{miplib2017}, and
Cor@l~\cite{linderoth2005noncommercial}.
These test sets consist of 1846 MINLP instances and 666 MILP instances, respectively.
After structure detection experiments, only those instances were chosen for performance evaluations that
either contain bilinear \firstrevision{terms} in the problem formulation, or where our algorithm derived
bilinear \firstrevision{terms}.
This resulted in
test sets of 1357 MINLP instances and 195 MILP instances.

The algorithms were implemented in the MINLP solver SCIP~\cite{BestuzhevaBesanconEtal2023_Enabling}.
We used a development branch of SCIP
(githash \texttt{dd6c54a9d7}) compiled with SoPlex~5.0.2.4 for LP solving,
CppAD 20180000.0 for automatic differentiation, 
PaPILO 1.0.0.1 for additional MILP presolving,
bliss 0.73p for graph automorphism computations for symmetry detection,
and Ipopt 3.12.11 for NLP solving.
The experiments were carried out on a cluster of Dell Poweredge M620 blades with 2.50GHz Intel Xeon
CPU E5-2670 v2 CPUs, with 2 CPUs and 64GB memory per node.
The time limit was set to one hour, the optimality gap tolerance to $10^{-4}$ for MINLP instances and
to $10^{-6}$ for MILP instances, and the following settings were used for all runs, where applicable:

\begin{itemize}
\item The maximum number of unknown bilinear terms that a product of a row and a bound factor can
have in order to be used was set to 20.
Unknown bilinear terms are those terms $x_ix_j$ for which no $w_{ij}$ variable exists in the
problem, or its extended formulation which SCIP constructs for the purposes of creating an LP
relaxation of an MINLP.
\item RLT cut separation was called every 10
nodes of the branch-and-bound tree.
\item In every non-root node where separation was called, 1 round of separation was performed. In
the root node, 10 separation rounds were performed.
\item Unless specified otherwise, implicit product detection and projection filtering were enabled
and the new separation algorithm \firstrevision{based on row marking} described in Section~\ref{subsection:marking} was used.
\end{itemize}

In SCIP, before RLT cuts \firstrevision{are} added to the solver's cut storage, they \firstrevision{pass} an additional
cut selection procedure.
The procedure itself \firstrevision{is} analogous to SCIP's default cut selection procedure that is based on cut
efficacy and cut parallelism.
The purpose of this additional round of cut selection \firstrevision{is} to compare RLT cuts between themselves
before allowing them to compete against other cuts generated by various solver components.


\subsection{Impact of RLT Cuts}\label{sec:impact_of_RLT}

In this section we evaluate the performance impact of RLT cuts.
The following settings were used:
\begin{itemize}
	\item \emph{Off} - RLT cuts are disabled;
	\item \emph{ERLT} - RLT cuts are added for products that exist explicitly in the problem; \firstrevision{row marking and projection filtering are enabled}.
	\item \emph{IERLT} - RLT cuts are added for both implicit and explicit products; \firstrevision{row marking and projection filtering are enabled}.
\end{itemize}
The setting \emph{ERLT} was used for the MINLP test set only, since MILP instances contain no
explicitly defined bilinear \firstrevision{terms}.

We report overall
numbers of instances, numbers
of solved instances, shifted geometric means of the runtime (shift 1 second),
and the number of nodes in the branch-and-bound tree (shift
100 nodes), and relative differences between settings.
\firstrevision{Clean instances are the instances where the solver was successful for all of the settings.
We report our results on clean instances only since it is unclear how to compute the contribution of
failed instances to the mean time and number of nodes.}

Additionally, we report results on subsets of \firstrevision{clean} instances.
Affected instances are instances where a change of setting leads to a difference in the solving
process, indicated by a difference in the number of LP iterations.
\bracket{x}{timelim} denotes the subset of
instances which took the solver at least $x$ seconds to solve with at least one setting, and were
solved to optimality with at least one setting.
All-optimal is the
subset of instances which were solved to optimality with both settings.

\firstrevision{Since the tables below only include clean instances, we will briefly discuss the
failures that lead to the exclusion of an instance.
There were 3 MILP instances (including permuted instances) where SCIP failed with setting \emph{Off} and 1 instance where it failed with setting \emph{IERLT}.
Among MINLP instances (including permuted instances), there were 155 instances where the solver failed with setting \emph{Off}, 161 instance with setting \emph{ERLT}, and 158 instances with setting \emph{IERLT}.
Such failures occur due to numerical errors and memory limits, and RLT cuts do not have a considerable effect on solver stability.
}

\begin{table}[h]
\caption{Impact of RLT cuts: MILP instances}
\label{tbl:implicit_milp}
\smallscript
\begin{tabular*}{\textwidth}{@{}l@{\;\;\extracolsep{\fill}}rrrrrrrrrrrr@{}}
\toprule
&           & \multicolumn{3}{c}{\emph{Off}} & \multicolumn{3}{c}{\emph{IERLT}} & \multicolumn{2}{c}{\emph{IERLT}/\emph{Off}} \\
\cmidrule{3-5} \cmidrule{6-8} \cmidrule{9-10}
Subset                & instances &       solved &       time &        nodes & solved    & time    & nodes   &        time &        nodes \\
\midrule
\firstrevision{Clean}                   &       971 &          905 &       \textbf{45.2} &         1339 &    \textbf{909}    & 46.7    &  \textbf{1310}   &        1.03 &         0.98 \\
Affected              &       581 &          571 &       \textbf{48.8} &         1936 &    \textbf{575}    & 51.2    &  \textbf{1877}   &        1.05 &         0.97 \\
\cmidrule{1-10}
\bracket{0}{timelim}    &       915 &          905 &       \textbf{34.4} &         1127  &   \textbf{909}    & 35.6    &  \textbf{1104}   &        1.04 &         0.98 \\
\bracket{1}{timelim}    &       832 &          822 &       \textbf{47.2} &         1451  &   \textbf{826}    & 49.0    &  \textbf{1420}   &        1.04 &         0.98 \\
\bracket{10}{timelim}   &       590 &          580 &      \textbf{126.8} &         3604  &   \textbf{584}    & 133.9   &  \textbf{3495}   &        1.06 &         0.97 \\
\bracket{100}{timelim}  &       329 &          319 &       439.1 &        9121 &    \textbf{323}    & \textbf{430.7}   &  \textbf{8333}   &        0.98 &         0.91 \\
\bracket{1000}{timelim} &        96 &           88 &      1436.7 &       43060 &    \textbf{92}     & \textbf{1140.9}  &  \textbf{31104}  &        0.79 &         0.72 \\
\cmidrule{1-10}
All-optimal           &       899 &          899 &       \textbf{31.9} &         \textbf{1033} &    899    & 34.1    &  1053   &        1.07 &         1.02 \\
\bottomrule
\end{tabular*}
\end{table}

Table~\ref{tbl:implicit_milp} shows the impact of RLT cuts on MILP performance.
We observe a slight increase in \firstrevision{mean} time \firstrevision{on all clean instances} when RLT cuts are enabled, and a slight decrease in
number of nodes.
The difference \firstrevision{in the numbers of nodes} is more pronounced on `difficult' instances: a $9\%$ decrease in number of nodes on
subset \bracket{100}{timelim} and $28\%$ on subset \bracket{1000}{timelim}.
\firstrevision{Notably, on these subsets RLT cuts decrease the mean time, leading to a small decrease of $2\%$ on subset \bracket{100}{timelim} and a more substantial decrease of $21\%$ on subset \bracket{1000}{timelim}.}

For MINLPs, Table~\ref{tbl:implicit_minlp_existing} reports the impact of RLT cuts derived from explicitly defined bilinear
products, comparing the settings \emph{Off} and \emph{ERLT}.
A substantial decrease in running times and tree sizes is observed across all
subsets, with a $15\%$ decrease in the mean time and a $19\%$ decrease in the number of nodes on
all \firstrevision{clean} instances, and a $87\%$ decrease in the mean time and a $88\%$
decrease in the number of nodes on the subset
\bracket{1000}{timelim}.
123 more instances are solved with \emph{ERLT} than with \emph{Off}.

\begin{table}[h]
\caption{Impact of RLT cuts derived from explicit products: MINLP instances}
\label{tbl:implicit_minlp_existing}
\smallscript
\begin{tabular*}{\textwidth}{@{}l@{\;\;\extracolsep{\fill}}rrrrrrrrrrrr@{}}
\toprule
&           & \multicolumn{3}{c}{\emph{Off}} & \multicolumn{3}{c}{\emph{ERLT}} & \multicolumn{2}{c}{\emph{ERLT}/\emph{Off}} \\
\cmidrule{3-5} \cmidrule{6-8} \cmidrule{9-10}
Subset                & instances &       solved &       time &        nodes & solved    & time    & nodes   &        time &        nodes \\
\midrule
\firstrevision{Clean}                   & 6622      &      4434    &   67.5     &      3375    &   \textbf{4557}    &  \textbf{57.5}   &  \textbf{2719}   &    0.85     &     0.81     \\
Affected              & 2018      &      1884    &   18.5     &      1534    &   \textbf{2007}    &  \textbf{10.6}   &  \textbf{784}   &    0.57     &     0.51     \\
\cmidrule{1-10}
\bracket{0}{timelim}  & 4568      &      4434    &   10.5     &      778     &   \textbf{4557}    &  \textbf{8.2}    &  \textbf{569}    &    0.78     &     0.73     \\
\bracket{1}{timelim}  & 3124      &      2990    &   28.3     &      2081    &   \textbf{3113}    &  \textbf{20.0}   &  \textbf{1383}   &    0.71     &     0.67     \\
\bracket{10}{timelim} & 1871      &      1737    &   108.3    &      6729    &   \textbf{1860}    &  \textbf{63.6}   &  \textbf{3745}   &    0.59     &     0.56     \\
\bracket{100}{timelim}  & 861       &      727     &   519.7    &     35991    &   \textbf{850}     &  \textbf{196.1}  &  \textbf{12873}  &    0.38     &     0.36     \\
\bracket{1000}{timelim} & 284       &      150     &  2354.8    &    196466    &   \textbf{273}     &  \textbf{297.6}  &  \textbf{23541}  &    0.13     &     0.12     \\
\cmidrule{1-10}
All-optimal           & 4423      &      4423    &   8.6      &       627    &   4423    &  \textbf{7.5}    &  \textbf{518}    &    0.87     &     0.83     \\
\bottomrule
\end{tabular*}
\end{table}

Table~\ref{tbl:implicit_minlp_hidden} evaluates the impact of RLT cuts derived from implicit
bilinear \firstrevision{terms}, comparing the settings \emph{ERLT} and \emph{IERLT}.
Similarly to MILP instances, the mean time slightly increases and the mean number of nodes
slightly decreases when additional RLT cuts are enabled, but on MINLP instances, the increase
in the mean time persists across different instance subsets and is most pronounced ($9\%$) on
the subset \bracket{100}{timelim}, and the number of nodes increases by
$6-7\%$ on subsets \bracket{100}{timelim} and \bracket{1000}{timelim}.
\firstrevision{Thus, the addition of implicit RLT cuts is overall slightly detrimental to MINLP performance.}

\begin{table}[h]
\caption{Impact of RLT cuts derived from implicit products: MINLP instances}
\label{tbl:implicit_minlp_hidden}
\smallscript
\begin{tabular*}{\textwidth}{@{}l@{\;\;\extracolsep{\fill}}rrrrrrrrrrrr@{}}
\toprule
&           & \multicolumn{3}{c}{\emph{ERLT}} & \multicolumn{3}{c}{\emph{IERLT}} & \multicolumn{2}{c}{\emph{IERLT}/\emph{ERLT}} \\
\cmidrule{3-5} \cmidrule{6-8} \cmidrule{9-10}
Subset                & instances &       solved &       time &        nodes & solved    & time    & nodes   &        time &        nodes \\
\midrule
\firstrevision{Clean}                   & 6622      &     4565     &    \textbf{57.0}    &       2686   &    \textbf{4568}   &  57.4   &  \textbf{2638}   &     1.01    &    0.98      \\
Affected              & 1738      &     1702     &    \textbf{24.2}    &       1567   &    \textbf{1705}   &  24.8   &  \textbf{1494}   &     1.02    &    0.95      \\
\cmidrule{1-10}
\bracket{0}{timelim}  & 4601      &     4565     &    \textbf{8.5}     &       587    &    \textbf{4568}   &  8.6    &  \textbf{576}    &     1.01    &    0.98      \\
\bracket{1}{timelim}  & 3141      &     3105     &    \textbf{21.1}    &       1436   &    \textbf{3108}   &  21.4   &  \textbf{1398}   &     1.01    &    0.97      \\
\bracket{10}{timelim} & 1828      &     1792     &    \textbf{74.1}    &       4157   &    \textbf{1795}   &  75.4   &  \textbf{4012}   &     1.02    &    0.97      \\
\bracket{100}{timelim}  & 706       &     670      &    \textbf{359.9}   &       \textbf{22875}  &    \textbf{673}    &  390.4  &  24339  &     1.09    &    1.06      \\
\bracket{1000}{timelim} & 192       &     156      &    \textbf{1493.3}  &       \textbf{99996}  &    \textbf{159}    &  1544.7 &  107006 &     1.03    &    1.07      \\
\cmidrule{1-10}
All-optimal           & 4532      &     4532     &    \textbf{7.7}     &       540    &    4532   &  7.8    &  \textbf{529}    &     1.02    &    0.98      \\
\bottomrule
\end{tabular*}
\end{table}

Table~\ref{tbl:rootnode} reports numbers of instances for which
a change in the root node dual bound was observed, where the relative difference is quantified as
$\frac{\gamma_2-\gamma_1}{\gamma_1}$,
where $\gamma_1$ and $\gamma_2$ are root node dual bounds obtained
with the first and second settings, respectively.
The range of the change is specified in the column `Difference', and each column shows numbers of
instances for which one or the other setting provided a better dual bound, within
given range.

The results of comparisons \emph{Off}/\emph{IERLT} for MILP instances and \emph{Off}/\emph{ERLT}
for MINLP instances are \firstrevision{mostly} consistent with the effect of RLT cuts on performance observed in
Tables~\ref{tbl:implicit_milp} and~\ref{tbl:implicit_minlp_existing}.
Interestingly, \emph{IERLT} performs better than \emph{ERLT} in terms of root node dual bound
quality.
Thus, RLT cuts derived from implicit products in MINLP instances tend to improve root node
relaxations, \firstrevision{although this does not always translate into an improvement in overall solver performance.}

\begin{table}[h]
\caption{Root node dual bound differences}
\label{tbl:rootnode}
\smallscript
\begin{tabular*}{\textwidth}{@{}l@{\;\;\extracolsep{\fill}}ccc@{}}
\toprule
           & \multicolumn{1}{c}{MILP} & \multicolumn{2}{c}{MINLP} \\
\cmidrule{2-2} \cmidrule{3-4}
Difference                & \emph{Off} / \emph{IERLT} &      \emph{Off}  /     \emph{ERLT} & \emph{ERLT}    / \emph{IERLT} \\
\midrule
0.01-0.2                   & 54      /      \textbf{62}      &    224    /      \textbf{505}    &   379      /  \textbf{441}      \\
0.2-0.5                    & 2       /      \textbf{4}       &    23     /      \textbf{114}    &   44       /  \textbf{48}       \\
0.5-1.0                    & 0       /      \textbf{3}       &    40     /      \textbf{150}    &   19       /  \textbf{30}       \\
$>$1.0                     & 0       /      \textbf{2}       &    4      /      \textbf{182}    &   4        /  \textbf{23}       \\
\bottomrule
\end{tabular*}
\end{table}

\firstrevision{Appendix~\ref{appendix1} provides detailed tables showing per-instance results for
subsets of MILP and MINLP instances.
Despite the performance variability which, as expected, is observed in these tables, it is clear
that for MINLP instances, settings \emph{ERLT} and \emph{IERLT} perform better than \emph{Off} for
most instances.
The comparison between settings \emph{Off} and \emph{IERLT} for MILP instances and \emph{ERLT} and
\emph{IERLT} for MINLP instances (that is, the impact of cuts for implicit products) is not as clear.
There are, however, groups of similarly structured instances where one setting clearly outperforms
the other: for example, cuts for implicit products are consistently detrimental for MILP instance
group \emph{csched*}, mostly detrimental for MILP instance group \emph{lectsched*}, mostly
beneficial for MILP instance group \emph{tr12-30*}, and mostly beneficial for MINLP instance group
\emph{smallinvSNPr*}.}

\firstrevision{Although in Tables~\ref{tbl:implicit_milp} and \ref{tbl:implicit_minlp_hidden} the
addition of RLT cuts for implicit products leads to slowdowns on most instance subsets, there exist
problems that clearly benefit from these cuts, and the decision on
whether to employ them depends on the problem.
For heterogeneous instance sets, setting \emph{IERLT} is the best performing only for challenging
MILP instances, setting \emph{Off} is the best performing for remaining MILP instances, and setting
\emph{ERLT} is the best performing for MINLP instances.}

\firstrevision{In the following sections, the baseline for comparisons is the setting \emph{IERLT},
which is compared to settings obtained by disabling or enabling individual features so as to evaluate
the features' effects.}

\subsection{Separation}

This section assesses the performance impact of the RLT cut separation algorithm \firstrevision{based on row marking} proposed in
this paper.
\firstrevision{Table~\ref{tbl:separation_milp} presents results with the following settings:}
\begin{itemize}
\item \firstrevision{\emph{IERLT} is identical to the setting \emph{IERLT} used throughout the section,
with RLT cuts enabled for explicit and implicit products, and row marking and projection filtering enabled.}
\item \firstrevision{\emph{No-mark} is identical to \emph{IERLT}, except that it employs standard separation instead
of the row marking algorithm described in Section~\ref{section:separation}.}.
\end{itemize}
From Table~\ref{tbl:separation_milp}, we can see that
row marking reduces the running time by 63\% on MILP
instances, by 70\% on affected MILP instances, by 12\% on MINLP instances and by 22\% on affected MINLP instances.
The number of nodes increases when row marking is enabled because, due to the
decreased separation time, the solver can explore more nodes before reaching the time limit: this is confirmed
by the fact that on the subset All-optimal,
the number of nodes remains nearly unchanged.

\begin{table}[h]
\caption{Separation algorithm comparison}
\label{tbl:separation_milp}
\smallscript
\begin{tabular*}{\textwidth}{@{}l@{\;\;\extracolsep{\fill}}rrrrrrrrrrrr@{}}
\toprule
&           & \multicolumn{3}{c}{\emph{No-mark}} & \multicolumn{3}{c}{\firstrevision{\emph{IERLT}}} & \multicolumn{2}{c}{\firstrevision{\emph{IERLT}}/\emph{No-mark}} \\
\cmidrule{3-5} \cmidrule{6-8} \cmidrule{9-10}
subset                & instances &       solved &       time  &        nodes & solved    & time    & nodes   &        time &        nodes \\
\midrule
\multicolumn{10}{c}{MILP:} \\
\midrule
\firstrevision{Clean} & 949    &    780    &    124.0    &    \textbf{952}    &   \textbf{890}      &  \textbf{45.2}   &  1297   &    0.37   &   1.37  \\
Affected              & 728    &    612    &    156.6    &    \textbf{1118}   &   \textbf{722}      &  \textbf{46.4}   &  1467   &    0.30   &     1.31     \\
All-optimal           & 774    &    774    &    58.4     &      \textbf{823}    &   774      &  \textbf{21.2}   &  829    &    0.36     &     1.01     \\
\cmidrule{1-10}
\multicolumn{10}{c}{MINLP:} \\
\midrule
\firstrevision{Clean} &  6546  &    4491   &   64.5      &    \textbf{2317}     &  \textbf{4530} & \textbf{56.4} & 2589   &   0.88      &   1.12       \\
Affected              &  3031  &    2949   &   18.5      &    \textbf{1062}     &   \textbf{2988} & \textbf{14.3} & 1116  &   0.78      &   1.05       \\
All-optimal           &  4448  &    4448   &   9.1       &    \textbf{494}      &   4448         & \textbf{7.4}  & 502    &   0.81      &   1.02       \\
\bottomrule
\end{tabular*}
\end{table}

Table~\ref{tbl:separation_times} analyzes the percentage of time that RLT cut
separation takes out of overall running time, showing the arithmetic mean and maximum over all
instances, numbers of instances for which the percentage was within a given interval, and numbers
of failures.
Instances where the total running time was less than or equal to 0.1 seconds were disregarded,
since, due to rounding, the computed percentages could be misleading.
The average percentage is reduced from $54.2\%$ to $2.8\%$ for MILP instances and from $15.1\%$ to
$2.4\%$ for MINLP instances, and the maximum percentage is reduced from $99.6\%$ to $71.6\%$ for
MILP instances and from $96.6\%$ to $66.4\%$ for MINLP instances.
The numbers of instances where RLT cut separation took less than $5\%$ of solving time is
considerably larger with \firstrevision{\emph{IERLT}}, and the numbers of instances where the percentage
belongs to any one of the remaining intervals have decreased.
The numbers of failures are reduced with \firstrevision{\emph{IERLT}}, mainly due to avoiding failures that occur when the
solver runs out of memory.
Overall, the percentage of time taken up by RLT cut separation is reduced
considerably by applying row marking.

\begin{table}[h]
\caption{Separation times}
\label{tbl:separation_times}
\smallscript
\begin{tabular*}{\textwidth}{@{}l@{\;\;\extracolsep{\fill}}lrrrrrrrrrrrr@{}}
\toprule
Test set &  Setting               & avg \% &       max \% &       N($<5\%$)  &   N(5-20\%)   & N(20-50\%)    & N(50-100\%) & fail \\
\midrule
MILP & \emph{No-mark}            &  54.6     &     99.6     &   114        &    117        &   166      & 550   &  16      \\
& \firstrevision{\emph{IERLT}}                  &  \textbf{2.8}     &     \textbf{71.6}     &   832      &    85      &   31    & 4  &  0      \\
\cmidrule{1-9}
MINLP &   \emph{No-mark}            &  16.7     &     96.6     &   2973      &    1237      &   1079    & 680  &  77      \\
& \firstrevision{\emph{IERLT}}                  &  \textbf{2.6}     &     \textbf{66.4}     &   5412      &    372      &   201    & 47  &  16      \\
\bottomrule
\end{tabular*}
\end{table}

\subsection{Redundant Product Filtering}

This section evaluates the impact of redundancy filtering introduced in
Section~\ref{sec:redchecks}.
Tables~\ref{tbl:redundancy_milp} and~\ref{tbl:redundancy_minlp}
compare performance between the following two settings:
\begin{itemize}
\item \firstrevision{\emph{IERLT}: identical} to the setting \firstrevision{\emph{IERLT}} used throughout the section,
\firstrevision{with RLT cuts enabled for explicit and implicit products, and row
marking and projection filtering enabled.}
In particular, the detection algorithm enforces condition (2) from
Section~\ref{subsection:detection_theory} and does not perform product relation redundancy checks.
\item \emph{RPfilter}: the detection algorithm does not enforce condition (2) from
Section~\ref{subsection:detection_theory}, but performs product relation redundancy checks
described in Section~\ref{sec:redchecks}.
The product is considered to be sufficiently non-redundant if
$\frac{\min\{\overline{x}^*_j,\overline{x}_j\} - \max\{\underline{x}^*_j,\underline{x}_j\}}{\overline{x}_j - \underline{x}_j} \geq 0.3$,
where $\underline{x}^*_j$ and $\overline{x}^*_j$ are, respectively, the smallest and the largest
values of $x_j$ such that the product relation dominates both linear constraints that it is derived
from.
\end{itemize}

Since the product filtering strategy affects the results of product detection, we ran detection
anew with setting \emph{RPfilter}.
To avoid unnecessary computation, we excluded some instances where the solver always reached the
time limit with all settings from the test sets.


\begin{table}
\caption{Product filtering comparison: MILP instances}
\label{tbl:redundancy_milp}
\smallscript
\begin{tabular*}{\textwidth}{@{}l@{\;\;\extracolsep{\fill}}rrrrrrrrrrrr@{}}
\toprule
&           & \multicolumn{3}{c}{\firstrevision{\emph{IERLT}}} & \multicolumn{3}{c}{\emph{RPfilter}} & \multicolumn{2}{c}{relative} \\
\cmidrule{3-5} \cmidrule{6-8} \cmidrule{9-10}
Subset                & instances &       solved  &   time     &   nodes          & solved         & time             & nodes           &  time   &   nodes  \\
\midrule
\firstrevision{Clean}                     &  1638     &    1004     &   560.1    &   \textbf{5563}  &  \textbf{1007} &  \textbf{551.6}  & 5813            &   0.99  &   1.05   \\
Affected                &  690      &    664      &   238.7    &   \textbf{4769}  &   \textbf{667} & \textbf{231.9}   & 4929            &   0.97  &   1.03   \\
\cmidrule{1-10}
\bracket{0}{timelim}    &  1030     &    1004     &   186.3    &   \textbf{2417}  &  \textbf{1007} & \textbf{181.8}   & 2472            &   0.98  &   1.02   \\
\bracket{1}{timelim}    &  1026     &    1000     &   189.7    &   \textbf{2440}  &  \textbf{1003} & \textbf{185.1}   & 2496            &   0.98  &   1.02   \\
\bracket{10}{timelim}   &  946      &    920      &   255.2    &   \textbf{3009}  &  \textbf{923}  & \textbf{248.2}   & 3083            &   0.97  &   1.02   \\
\bracket{100}{timelim}  &  683      &    657      &   537.0    &   \textbf{5251}  &  \textbf{660}  & \textbf{519.7}   & 5306            &   0.97  &   1.01   \\
\bracket{1000}{timelim} &  248      &    222      &   1657.0   &   24108          &   \textbf{225} & \textbf{1621.2}  &  \textbf{23967} &   0.98  &   0.99   \\
\cmidrule{1-10}
All-optimal           &   981       &     981     &   165.3    &   \textbf{2049}  &   981          & \textbf{161.8}   &  2125           &   0.98  &   1.04   \\
\bottomrule
\end{tabular*}
\end{table}

\begin{table}
\caption{Product filtering comparison: MINLP instances}
\label{tbl:redundancy_minlp}
\smallscript
\begin{tabular*}{\textwidth}{@{}l@{\;\;\extracolsep{\fill}}rrrrrrrrrrrr@{}}
\toprule
&           & \multicolumn{3}{c}{\firstrevision{\emph{IERLT}}} & \multicolumn{3}{c}{\emph{RPfilter}} & \multicolumn{2}{c}{relative} \\
\cmidrule{3-5} \cmidrule{6-8} \cmidrule{9-10}
Subset                  & instances &  solved       &   time           &  nodes          & solved & time            & nodes          & time  &  nodes \\
\midrule
\firstrevision{Clean}                     &  708      &  \textbf{598} &  52.9            &  \textbf{2412}  & 592    &  \textbf{51.9}  &  2413          & 0.98  &  1.00  \\
Affected                &  288      &  \textbf{586} &  44.8            &  \textbf{2019}  & 580    &  \textbf{43.0}  &  2035          & 0.96  &  1.01  \\
\cmidrule{1-10}
\bracket{0}{timelim}    &  600      &  \textbf{598} &  24.3            &  \textbf{1503}  & 592    &  \textbf{23.7}  &  1509          & 0.98  &  1.00  \\
\bracket{1}{timelim}    &  566      &  \textbf{564} &  28.8            &  \textbf{1759}  & 558    &  \textbf{28.1}  &  1766          & 0.98  &  1.00  \\
\bracket{10}{timelim}   &  352      &  \textbf{350} &  96.3            &  \textbf{4698}  & 344    &  \textbf{92.4}  &  4749          & 0.96  &  1.01  \\
\bracket{100}{timelim}  &  163      &  \textbf{161} &  412.8           &  \textbf{14834} &  155   &  \textbf{381.8} &  14912         & 0.93  &  1.01  \\
\bracket{1000}{timelim} &  50       &  \textbf{48}  &  \textbf{1342.5} &  \textbf{71128} &  42    &  1477.8         &  78552         & 1.10  &  1.10  \\
\cmidrule{1-10}
All-optimal             &  590      &  590          &  22.7            &  1393           &  590   &  \textbf{22.0}  &  \textbf{1387} & 0.97  &  1.00  \\
\bottomrule
\end{tabular*}
\end{table}

We observe a consistent decrease of the running time on both MILP and MINLP instances on almost all
instance subsets.
On MILP instances, the new product filtering strategy yields a running time decrease of $3\%$ on
affected instances, and on MINLP instances, the decrease in running time is $4\%$ on affected
instances and $7\%$ on the subset \bracket{100}{timelim}, where it has the largest impact on
running time.
The only exception is the subset \bracket{1000}{timelim} of 50 MINLP instances, where setting
\emph{RPfilter} causes a slowdown of $10\%$.
The relative difference in the number of nodes tends to be smaller than the relative difference in
running time, with a slight increase occurring on most instance subsets.
3 more MILP instances and 4 less MINLP instances are solved with setting \emph{RPfilter}.

Additionally, we tried storing products that failed the nonredundancy check and using them for
substitution instead of relaxations in the linearization step, but ignoring them in row marking.
This strategy was nearly performance neutral compared to the setting \firstrevision{\emph{IERLT}}.

\subsection{Projection}

In this section we demonstrate the impact of projection filtering  on solver performance.
Tables~\ref{tbl:projection_milp} and \ref{tbl:projection_minlp} compare performance with \firstrevision{the following settings}:
\begin{itemize}
\item \firstrevision{\emph{IERLT}: identical to the setting \emph{IERLT} used throughout the section,
with RLT cuts enabled for explicit and implicit products, and row marking and projection filtering enabled.}
\item \firstrevision{\emph{No-proj}: identical to IERLT, except projection filtering described in
Section~\ref{subsection:projection} is disabled.}
\end{itemize}

Projection filtering has a minor impact on performance.
When comparing the runs where projection filtering is disabled and enabled, the relative difference
in time and nodes does not exceed $1\%$ on both MILP and MINLP instances, with time slightly
decreasing when projection filtering is applied on all subsets except for \bracket{1000}{timelim},
where it increases by 2 and 3\%, respectively.
The difference in number of nodes is more pronounced on
except for affected MILP instances where
projection filtering decreases the number of nodes by 4\%.
This is possibly occurring due to the
effect of prioritizing McCormick inequalities to RLT cuts when enforcing derived product relations.
The number of solved instances remains almost unchanged, with one less instance being solved on
both MILP and MINLP test sets when projection filtering is enabled.

On MINLP instances, the impact of projection filtering is smaller; it improves performance by
1-3\% on all subsets except for the most challenging MINLP instances, where it causes a slowdown of
7\%.
We analyzed the difference in time taken up by separation, but we do not report it in a table
since it is negligibly small.
The difference in the average percentage of total time
spent in RLT cut separation did not exceed 0.3 percentage points for both MILP and MINLP instances.

\begin{table}
\caption{Projection impact evaluation: MILP instances}
\label{tbl:projection_milp}
\smallscript
\begin{tabular*}{\textwidth}{@{}l@{\;\;\extracolsep{\fill}}lrrrrrrrrrrrr@{}}
\toprule
&      & \multicolumn{3}{c}{\emph{No-proj}} & \multicolumn{3}{c}{\firstrevision{\emph{IERLT}}} & \multicolumn{2}{c}{\firstrevision{\emph{IERLT}}/\emph{No-proj}} \\
\cmidrule{3-5} \cmidrule{6-8} \cmidrule{9-10}
subset                & instances &       solved &       time  &        nodes & solved    & time    & nodes   &        time &        nodes \\
\midrule
\firstrevision{Clean}                   &   972    &     \textbf{912}       &    46.4     &    1329      &   911     &  \textbf{46.1}   &  \textbf{1302}   &   0.99      &     0.98     \\
Affected              &   530    &     \textbf{523}       &    75.7     &    3092      &   522     &  \textbf{74.6}   &  \textbf{2964}   &   0.99      &     0.96     \\
\cmidrule{1-10}
\bracket{0}{timelim}  &   919    &     \textbf{912}       &    36.0     &    1155      &   911     &  \textbf{35.7}   &  \textbf{1126}   &   0.99      &     0.98     \\
\bracket{1}{timelim}  &   832    &     \textbf{825}       &    50.3     &    1504      &   824     &  \textbf{49.8}   &  \textbf{1462}   &   0.99      &     0.97     \\
\bracket{10}{timelim} &   582    &     \textbf{575}       &    143.4    &    3886      &   574     &  \textbf{141.7}  &  \textbf{3741}   &   0.99      &     0.96     \\
\bracket{100}{timelim}  &   323    &     \textbf{316}       &    485.0    &    9601      &   315     &  \textbf{471.3}  &  \textbf{9065}   &   0.97      &     0.94     \\
\bracket{1000}{timelim} &    96    &      \textbf{89}       &    \textbf{1483.8}   &    45276     &   88      &  1512.2 &  \textbf{43061}  &   1.02      &     0.95     \\
\cmidrule{1-10}
All-optimal           &   904    &     904       &    33.5     &    1054      &   904     &  \textbf{33.4}   &  \textbf{1040}   &   1.00      &     0.99     \\
\bottomrule
\end{tabular*}
\end{table}

\begin{table}
\caption{Projection impact evaluation: MINLP instances}
\label{tbl:projection_minlp}
\smallscript
\begin{tabular*}{\textwidth}{@{}l@{\;\;\extracolsep{\fill}}rrrrrrrrrrrr@{}}
\toprule
&           & \multicolumn{3}{c}{\emph{No-proj}} & \multicolumn{3}{c}{\firstrevision{\emph{IERLT}}} & \multicolumn{2}{c}{\firstrevision{\emph{IERLT}}/\emph{No-proj}} \\
\cmidrule{3-5} \cmidrule{6-8} \cmidrule{9-10}
Subset                & instances &       solved &       time  &        nodes & solved    & time    & nodes   &        time &        nodes \\
\midrule
\firstrevision{Clean}                   &   6637    &     \textbf{4582}     &    57.9     &    2689      &  4581     &  \textbf{57.7}   &  \textbf{2674}   &   1.00      &    0.99      \\
Affected              &   2476    &     \textbf{2438}     &    23.3     &    1681      &  2437     &  \textbf{23.1}   &  \textbf{1660}   &   0.99      &    0.99      \\
\cmidrule{1-10}
\bracket{0}{timelim}  &   4620    &     \textbf{4582}     &    8.8      &    595       &  4581     &  \textbf{8.7}    &  \textbf{590}    &   0.99      &    0.99      \\
\bracket{1}{timelim}  &   3137    &     \textbf{3099}     &    22.4     &    1483      &  3098     &  \textbf{22.3}   &  \textbf{1467}   &   0.99      &    0.99      \\
\bracket{10}{timelim} &   1854    &     \textbf{1816}     &    77.7     &    \textbf{4253}      &  1815     &  \textbf{76.4}   &  4256   &   1.00      &    0.98      \\
\bracket{100}{timelim}  &   743     &     \textbf{705}      &    377.4    &    23389     &  704      &  \textbf{364.4}  &  \textbf{22680}  &   0.97      &    0.97      \\
\bracket{1000}{timelim} &   205     &     \textbf{167}      &    \textbf{1434.5}   &   \textbf{98443}     &  166      &  1480.7 &  105546 &   1.03      &    1.07      \\
\cmidrule{1-10}
All-optimal           &   4543    &     4543     &    8.0      &    539       &  4543     &  \textbf{7.9}    &  \textbf{533}    &   0.99      &    0.99      \\
\bottomrule
\end{tabular*}
\end{table}

In Section~\ref{subsection:projection}, we discussed the relation between projection filtering and
the strategy for adding McCormick envelopes: if McCormick cuts are added dynamically, then
projection filtering may result in some violated RLT cuts not being generated.
In order to check whether this affects the usefulness of projection filtering, we conducted a
control experiment on a reduced set of instances, forcing the addition of all McCormick inequalities to the LP relaxation.
Our results show that projection filtering remains nearly performance neutral, with an increase in running
time of $2\%$ on all instances and an increase in nodes of $1\%$ on all instances in the test set,
and a slightly more pronounced increase in running time and nodes on difficult instances.

Thus, projection filtering does not become more beneficial in the presence of all McCormick
inequalities.
This also indicates that projection filtering may prevent useful cuts from being added even in this
setup.
Since projection filtering only leads to ignoring cuts that would not be violated by the
current LP solution and thus would not be added at that node, the only impact it may have on the
solution process stems from some potentially useful cuts not being added to the cut pool.

\subsection{Performance Impact by Instance Size}

Table~\ref{tbl:group_times} analyzes the performance impact of features relating to RLT cuts
on both MILP and MINLP instances grouped by instance size, where size is defined as the number of
nonzeroes in the instance, using the geometric mean of solving time as performance measure.
We compare the setting \firstrevision{\emph{IERLT}}, where implicit products, marking and projection are enabled, to
settings \emph{Off}, \emph{No-mark} and \emph{No-proj} as defined in previous sections.

The impact of all features is not monotonic in instance size.
All \firstrevision{considered} features together (\firstrevision{\emph{IERLT}}) have the most impact on small-to-medium sized instances
(mean time reduced by 20\% compared to \emph{Off} on instances with 10-100 nonzeroes, and by 25\%
on instances with 100-1000 nonzeroes).
Row marking has the most impact (32-49\%) on instances with over 1000 nonzeroes.
Projection filtering has little impact on most subsets, but, notably, it yields an improvement of 10\% on the
largest instances (over 100000 nonzeroes).


\begin{table}
\caption{Geometric Mean Running Time for Instance Groups}
\label{tbl:group_times}
\smallscript
\begin{tabular*}{\textwidth}{@{}l@{\;\;\extracolsep{\fill}}rrrrrrr@{}}
\toprule
Setting / size              &  0-10  &  10-100        &  100-1000       &  1000-10000      &  10000-100000    &  100000-1000000 \\
\midrule
\emph{Off}                  &  1.3   &  3             &  47.3	        &  162.6           &  451.4           &	 1532.8 \\
\emph{IERLT}                &  1.3   &  \textbf{2.4}  &  35.7           &  \textbf{150.1}  &  \textbf{450.5}  &  \textbf{1431.7}     \\
\emph{No-mark}              &  1.3   &  2.8           &  38.4	        &  233.6           &  883.8           &	 2090.2	 \\
\emph{No-proj}              &  1.3   &  \textbf{2.4}  &  \textbf{35.0}  &  152.5           &  463.5           &	 1591.8	 \\
\emph{IERLT}/\emph{Off}     &  1     &  0.8           &  0.75	        &  0.92            &  1.00            &	 0.93	 \\
\emph{IERLT}/\emph{No-mark} &  1     &  0.86          &  0.93	        &  0.64            &  0.51            &	 0.68	 \\
\emph{IERLT}/\emph{No-proj} &  1     &  1             &  1.02	        &  0.98            &  0.97            &	 0.90 \\
\bottomrule
\end{tabular*}
\end{table}

\subsection{Experiments with Gurobi}

In this section we present results obtained by running
the mixed-integer quadratically-constrained programming solver Gurobi 10.0 beta~\cite{gurobi}.
The algorithms for implicit product detection and RLT cut separation are the same as in SCIP, although
implementation details may differ between the solvers.

The internal Gurobi test set was used, comprised of models sent by Gurobi customers and models from
public benchmarks,
chosen in a
way that avoids overrepresenting any particular problem class.
Whenever RLT cuts were enabled, so was implicit product detection, row marking and projection
filtering.
The time limit was set to 10000 seconds.

Table~\ref{tbl:gurobi} shows, for both MILP and MINLP test sets, the numbers of instances in the
test sets and their subsets, and the ratios of shifted geometric means of running time and number of
nodes of the runs with RLT cuts enabled, to the same means obtained with RLT cuts disabled.
The last row shows the numbers of instances solved with one setting and unsolved with the
other, that is, for example, ``RLT off: +41'' means that 41 instances were solved with the setting
``off'' that were not solved with the setting ``on''.

While the results cannot be directly compared to those obtained with SCIP due to the differences in
the experimental setup, we observe the same tendencies.
In particular, RLT cuts yield small improvements on MILP instances which become more pronounced on
subsets \bracket{100}{timelim} and \bracket{1000}{timelim}, and larger improvements are observed
on MINLP instances both in terms of geometric means and numbers of solved instances.
Relative differences are comparable to those observed with SCIP, but the impact of RLT cuts is
larger in Gurobi, and no slowdown is observed with Gurobi on any subset of MILP instances.


\begin{table}[h]
\caption{Results obtained with Gurobi 10.0 beta}
\label{tbl:gurobi}
\smallscript
\begin{tabular*}{\textwidth}{@{}l@{\;\;\extracolsep{\fill}}rrrrrrrrr@{}}
\toprule
&            \multicolumn{3}{c}{MILP} & \multicolumn{3}{c}{MINLP} \\
\cmidrule{2-4} \cmidrule{5-7}
Subset                & instances  &       timeR   &       nodeR     &    instances   &    timeR &        nodeR \\
\midrule
All                     & 5011     &      0.99     &     0.97  &        806         &   0.73   &  0.57   \\
\bracket{0}{timelim}    & 4830     &      0.99     &     0.96  &        505         &   0.57   &  0.44  \\
\bracket{1}{timelim}    & 3332     &      0.98     &     0.96  &        280         &   0.40   &  0.29  \\
\bracket{10}{timelim}   & 2410     &      0.97     &     0.93  &        188         &   0.29   &  0.20  \\
\bracket{100}{timelim}  & 1391     &      0.95     &     0.91  &        114         &   0.17   &  0.11  \\
\bracket{1000}{timelim} & 512      &      0.89     &     0.83  &        79          &   0.12   &  0.08  \\
\midrule
Solved               & \multicolumn{3}{c}{RLT off: +41; RLT on: +37} & \multicolumn{3}{c}{RLT off: +2; RLT on: +35} \\
\bottomrule
\end{tabular*}
\end{table}

Projection filtering had a negligibly small impact on performance on the MILP and MINLP test sets.





\section{Conclusion}\label{section:conclusion}

We developed a new RLT cut generation technique for MILPs, based on recognizing hidden nonlinear
structures represented by pairs of MILP constraints with binary variables.
We enhanced this detection by
filtering strategies that remove redundant and nearly redundant product relations, and designed a new RLT
cutting plane separation algorithm for MILPs and MINLPs.
We tested these methods with independent implementations in two solvers: the MICQP solver Gurobi, and the MINLP solver SCIP.

Our extensive computational study shows that
RLT cuts yield a considerable performance improvement for MINLP problems
and a small performance improvement for MILP problems which becomes more pronounced
for challenging instances.
The new separation algorithm \firstrevision{based on marking of promising rows} drastically reduces the computational burden of
RLT cut separation and is essential to an efficient implementation of RLT
cuts, enabling the speedups we observed when activating RLT.

\backmatter

\section*{Declarations}

\paragraph{Funding} The work for this article has been conducted within the Research Campus MODAL funded by the German
Federal Ministry of Education and Research (BMBF grant numbers 05M14ZAM, 05M20ZBM).

\paragraph{Competing interests} The authors confirm that they have no competing interests that are directly or indirectly related to this work.








\bibliography{rlt}

\newpage

\begin{appendices}
\firstrevisioncolor{
\section{Detailed Tables}\label{appendix1}

Tables~\ref{tbl:detailed_milp} and \ref{tbl:detailed_minlp} show detailed per-instance
comparisons of settings \emph{Off}, \emph{ERLT} and \emph{IERLT} on subsets of MILP and MINLP
instances, respectively.

\begin{smallscript}
\begin{longtable}{@{}p{8em}@{\;\;\extracolsep{\fill}}|rr|rr@{}}
\caption{Time and nodes on a subset of MILP instances with 4 permutations where at least one setting solved the instance, and at least one setting took $>1000$s to solve. For instances where the solver failed, time and node are displayed as $-1$.}
\label{tbl:detailed_milp} \\
                   & \multicolumn{2}{c}{Time} & \multicolumn{2}{c}{Nodes} \\
Instance           & \emph{Off}  &  \emph{IERLT}  &  \emph{Off}  &  \emph{IERLT} \\
\midrule
\endfirsthead
\caption{Time and nodes on a subset of MILP instances: continued} \\
                   & \multicolumn{2}{c}{Time} & \multicolumn{2}{c}{Nodes} \\
Instance           & \emph{Off}  &  \emph{IERLT}  &  \emph{Off}  &  \emph{IERLT} \\
\midrule
\endhead
\midrule
\endfoot
\endlastfoot
 csched008         &  \bf{520.8} & 1175.1 &   \bf{40452} &   87569 \\
 csched008.1       &  \bf{583.1} & 1145.9 &   \bf{39058} &  111636 \\
 csched008.3       &  \bf{967.8} & 1127.7 &   \bf{72877} &   83449 \\
 csched008.4       & \bf{1427.7} & 1544.3 &   \bf{93209} &  105273 \\
 danoint.2         & 3600.0 & \bf{3534.6} &  \bf{717873} &  783207 \\
 fiball.4          &  \bf{595.2} & 1390.3 &   \bf{13067} &   51828 \\
 glass4            & 1015.1 &  \bf{494.9} & 1013338 &  \bf{392765} \\
 glass4.2          & 1050.9 &  \bf{132.7} &  819774 &   \bf{68352} \\
 glass4.3          &  \bf{502.4} & 1967.5 &  \bf{313885} & 1100570 \\
 glass4.4          &  \bf{495.7} & 3600.0 &  \bf{370323} & 2889692 \\
 lectsched-1.1     & 1244.7 &  \bf{768.5} &  157685 &   \bf{21321} \\
 lectsched-1.3     &  \bf{256.3} & 3600.0 &   \bf{15759} &  273098 \\
 lectsched-3       &  \bf{139.4} & 1444.5 &     \bf{552} &    3276 \\
 lectsched-3.3     &  \bf{475.7} & 1471.7 &    2805 &    \bf{1825} \\
 map06             & \bf{1112.8} & 1157.8 &     \bf{839} &     918 \\
 map06.2           & 1192.6 & \bf{1034.0} &    1090 &     \bf{620} \\
 map06.4           & 1000.2 &  \bf{866.0} &     820 &     \bf{726} \\
 map10.1           & 1360.3 & \bf{1117.9} &    2154 &    \bf{1212} \\
 map10.3           & \bf{1078.4} & 1230.7 &    1804 &    \bf{1479} \\
 map10.4           &  \bf{692.0} & 1007.1 &     \bf{989} &    1507 \\
 neos-1053234.4    &  \bf{851.6} & 2426.5 &  \bf{192430} &  572908 \\
 neos-1171692.2    & 3600.0 &   \bf{10.1} &    5282 &       \bf{8} \\
 neos-1215891.2    & 1435.2 & \bf{1334.2} &    6136 &    \bf{4928} \\
 neos-555424       & 1318.0 &  \bf{646.7} &  101945 &   \bf{29371} \\
 neos-555424.1     & 3600.0 & \bf{1383.6} &  241941 &   \bf{79651} \\
 neos-555927.3     &  \bf{365.2} & 2431.2 &   \bf{41830} &  397695 \\
 neos-603073       & 3600.0 &  \bf{111.9} & 1096675 &   \bf{14655} \\
 neos-691058.2     & 2621.7 & \bf{1223.8} &    3618 &     \bf{628} \\
 neos-691073.2     & \bf{1318.9} & 1806.7 &    1180 &    \bf{1050} \\
 neos-693347       & 2272.8 & \bf{1930.7} &   17952 &   \bf{10093} \\
 neos-693347.1     &  \bf{978.5} & 2254.9 &   \bf{11526} &   19765 \\
 neos-693347.2     & \bf{1040.3} & 1499.0 &    \bf{8805} &   10551 \\
 neos-693347.3     & \bf{1341.7} & 2115.3 &    \bf{5673} &   17285 \\
 neos-693347.4     & \bf{1943.6} & 2465.2 &   \bf{13973} &   27587 \\
 neos-847051       & 3600.0 &    \bf{9.3} &  555907 &     \bf{549} \\
 neos-847051.1     & -1     &    \bf{6.0} &      -1 &       \bf{4} \\
 neos-847051.2     & 3600.0 &    \bf{5.2} &  985289 &       \bf{6} \\
 neos-847051.3     & -1     &    \bf{3.0} &  -1     &       \bf{3} \\
 neos-847051.4     & 3600.0 &    \bf{4.1} &  931598 &      \bf{14} \\
 neos-885086       & \bf{3382.8} & 3600.0 &     \bf{244} &     791 \\
 neos-885086.1     & \bf{2102.4} & 3600.0 &     \bf{201} &     432 \\
 neos-885086.2     & 2569.0 &   \bf{43.4} &     296 &       \bf{1} \\
 neos-911880       &  \bf{245.2} & 2421.4 &  \bf{109296} & 1153279 \\
 neos-911880.1     & 1388.0 &  \bf{290.6} &  525517 &   \bf{56546} \\
 neos-911880.2     &  \bf{253.9} & 3032.9 &   \bf{71993} &  423839 \\
 neos-911880.4     &  \bf{277.2} & 1346.7 &   \bf{56488} & 1038415 \\
 neos-911970.1     & 3600.0 &  \bf{441.4} & 2740634 &  \bf{244017} \\
 neos-911970.2     & 3600.0 &   \bf{64.4} &  883807 &    \bf{9063} \\
 neos-911970.4     & 3213.5 &  \bf{351.7} & 2801567 &   \bf{77945} \\
 neos-916792.3     & \bf{1173.1} & 1495.9 &   \bf{94156} &  125717 \\
 neos-916792.4     & \bf{1005.4} & 1014.9 &   \bf{85826} &   96293 \\
 net12.3           & 1243.4 &  \bf{432.1} &    2590 &     \bf{982} \\
 ns894788          & 2988.8 & \bf{1995.4} &  199994 &  \bf{102238} \\
 ns894788.2        & 1159.1 &  \bf{789.0} &   53436 &   \bf{39058} \\
 ns894788.3        & \bf{2546.2} & 3600.0 &  \bf{128593} &  173212 \\
 ns894788.4        &  \bf{461.3} & 3600.0 &   \bf{20684} &  188662 \\
 ran14x18-disj-8   & 2051.5 & \bf{1559.6} &  496355 &  \bf{271683} \\
 ran14x18-disj-8.1 & \bf{1298.8} & 1813.2 &  \bf{200082} &  335122 \\
 ran14x18-disj-8.2 & 2291.9 & \bf{1796.9} &  366296 &  \bf{355104} \\
 ran14x18-disj-8.3 & \bf{1288.4} & 1731.1 &  \bf{157630} &  318971 \\
 ran14x18-disj-8.4 & 1802.7 & \bf{1398.6} &  352915 &  \bf{170365} \\
 tr12-30           & 1461.1 &  \bf{804.0} &  650318 &  \bf{326771} \\
 tr12-30.1         & 1164.5 &  \bf{933.3} &  490969 &  \bf{381365} \\
 tr12-30.2         & 1044.0 &  \bf{863.9} &  443383 &  \bf{335560} \\
 tr12-30.3         &  \bf{783.1} & 1119.0 &  \bf{326205} &  474021 \\
 uct-subprob       & 3600.0 & \bf{2285.9} &  118040 &   \bf{65176} \\
 uct-subprob.1     & \bf{3378.5} & 3423.3 &  100452 &   \bf{92383} \\
 uct-subprob.2     & 2680.7 & \bf{2421.2} &   74647 &   \bf{65739} \\
 uct-subprob.3     & \bf{1713.0} & 2347.1 &   \bf{49204} &   68202 \\
 uct-subprob.4     & 2650.3 & \bf{1932.2} &   66351 &   \bf{56593} \\
 umts              & \bf{1405.0} & 1618.3 &  \bf{191410} &  229720 \\
 umts.1            & \bf{2175.3} & 3034.5 &  \bf{290634} &  455735 \\
 umts.2            & \bf{1507.5} & 3110.2 &  \bf{202894} &  462528 \\
 umts.3            & 1753.9 & \bf{1650.8} &  \bf{247547} &  256442 \\
 umts.4            & \bf{1761.7} & 2545.2 &  \bf{281500} &  448767 \\
\bottomrule
\end{longtable}
\end{smallscript}
}

\begin{table}[h!]
\caption{Time and nodes on a subset of MINLP instances with no permutations where at least one setting solved the instance, and at least one setting took $>1000$s to solve. For instances where the solver failed, time and node are displayed as $-1$. The full names of the first two instances are 0180404\_1d15min\_RmLinDisChar\_Feasib and AP\_20180404\_1d30min\_RmLinDis\_Feasib.}
\label{tbl:detailed_minlp}
\firstrevisioncolor
\smallscript
\begin{tabular*}{\textwidth}{@{}l@{\;\;\extracolsep{\fill}}|rrr|rrr@{}}
\toprule
                                        &            \multicolumn{3}{c}{Time} & \multicolumn{3}{c}{Nodes} \\
Instance                                & \emph{Off}  &  \emph{ERLT}  &  \emph{IERLT}  &  \emph{Off}  &  \emph{ERLT}  &  \emph{IERLT} \\
\midrule
0180404\_1d15min\_RmLinDisChar  & 2118.5 & \bf{1294.8} & 2115.6 &   18582 &   \bf{11885} &   28959 \\
AP\_20180404\_1d30min\_RmLinDis &  749.2 &  \bf{725.1} & 3600.0 &    \bf{9374} &   13832 &   59347 \\
autocorr\_bern40-05                     & 2639.8 & 2674.3 &  \bf{987.2} & 1168978 & 1168978 &  \bf{399615} \\
batchs151208m                        & -1 & -1 &   \bf{20.4} &       -1 &       -1 &     \bf{696} \\
bayes2\_50                           & 3600.0 & \bf{2955.7} & 2976.0 &  140686 &   \bf{94961} &   \bf{94961} \\
chp\_shorttermplan2d                 & 3600.0 &   11.7 &   \bf{11.4} &    4600 &       \bf{1} &       \bf{1} \\
crudeoil\_pooling\_ct2                & 3600.0 &   24.7 &   \bf{23.4} & 1836988 &     577 &     \bf{423} \\
ex1252                              & 3600.0 & 3600.0 &   \bf{39.4} & 1976844 & 1546510 &   \bf{21280} \\
ex1252a                             & 3600.0 &   \bf{14.4} & 3600.0 & 3327696 &    \bf{8174} & 1500306 \\
ex5\_2\_5                             & 3600.0 &    \bf{1.5} &    \bf{1.5} & 1397984 &     \bf{161} &     \bf{161} \\
flay06m                             & 3331.0 & \bf{3328.1} & 3569.2 & 3154885 & 3154885 & \bf{3094105} \\
fo9\_ar2\_1                           &  \bf{785.7} &  788.8 & 2042.7 & \bf{1011164} & \bf{1011164} & 2525951 \\
forest                              & 3600.0 &  \bf{140.8} &  143.8 & 1432689 &   \bf{33875} &  116131 \\
gabriel02                           & 1448.2 &  \bf{584.4} & 1563.6 &  192366 &   \bf{60336} &  205237 \\
gastrans135                         & 3600.0 &   88.9 &   \bf{81.8} &    5491 &      72 &      \bf{47} \\
gastrans582\_cool14                  & 3600.0 &   \bf{19.6} &   62.1 &  368396 &     \bf{180} &    2872 \\
graphpart\_clique-50                 &  926.1 &  \bf{918.1} & 1028.2 &    \bf{4786} &    \bf{4786} &    5149 \\
heatexch\_spec1                      &   57.1 & 3600.0 &    \bf{1.1} &   33563 & 1825713 &       \bf{2} \\
heatexch\_spec3                      & \bf{1127.3} & 3600.0 & 3600.0 &  \bf{147347} & 1201993 &  325967 \\
kall\_circles\_c8a                    & \bf{1139.7} & 2691.1 & 2984.3 &  \bf{630146} & 1041754 & 1041754 \\
kall\_diffcircles\_10                 &  \bf{973.8} & 1392.1 & 1449.9 &  645789 &  \bf{552883} &  \bf{552883} \\
multiplants\_mtg1a                   & 2649.2 & \bf{2643.1} & 3600.0 &  \bf{395907} &  493545 &  495304 \\
o7\_ar4\_1                            & \bf{1141.7} & 1148.3 & 2081.2 & \bf{1221378} & \bf{1221378} & 2360981 \\
qspp\_0\_12\_0\_1\_10\_1                  & \bf{1045.8} & 1061.4 & 1107.7 &  209746 &  209746 &  \bf{208806} \\
rsyn0815m04h                        & \bf{2746.9} & 2787.5 & 3600.0 &  \bf{418255} &  \bf{418255} &  512474 \\
rsyn0820m04hfsg                     & 3314.1 & \bf{3269.6} & 3600.0 &  621462 &  621462 &  \bf{516930} \\
rsyn0840h                           & 3600.0 & 3600.0 &  \bf{777.6} &   68992 &   68994 &   \bf{67358} \\
smallinvSNPr1b100-110               & 3600.0 &   89.1 &   \bf{88.2} &   55332 &     \bf{661} &     \bf{661} \\
smallinvSNPr1b150-165               & 3600.0 &  \bf{107.9} &  108.2 &   31107 &    \bf{2303} &    \bf{2303} \\
smallinvSNPr1b200-220               & 3600.0 &   68.7 &   \bf{67.1} &   22160 &     \bf{589} &     \bf{589} \\
smallinvSNPr2b100-110               & 2012.9 &   83.1 &   \bf{83.0} &   25464 &     \bf{637} &     \bf{637} \\
smallinvSNPr2b150-165               & 3600.0 &   79.9 &   \bf{77.5} &   36636 &    \bf{2600} &    \bf{2600} \\
smallinvSNPr2b200-220               & 3600.0 &   52.4 &   \bf{51.7} &   27038 &     \bf{698} &     \bf{698} \\
smallinvSNPr4b200-220               & 2885.3 &   \bf{46.8} &   47.6 &   30271 &     \bf{500} &     \bf{500} \\
sonet17v4                           & 3213.5 & \bf{3205.2} & 3235.8 &    4301 &    4301 &    \bf{3694} \\
sonet18v6                           & 3600.0 & 3600.0 & \bf{1639.6} &    3123 &    3124 &    \bf{2647} \\
syn20m04h                           & 2634.7 & \bf{2614.2} & 3600.0 &  \bf{690187} &  \bf{690187} &  927703 \\
syn30m04hfsg                        & 2289.3 & \bf{2267.9} & 3600.0 &  \bf{446351} &  \bf{446351} &  532670 \\
tln6                                & 1779.5 &    3.3 &    \bf{3.1} &  657078 &     \bf{915} &     \bf{915} \\
tln7                                & 3600.0 & \bf{1880.6} & 1890.0 &  \bf{338846} & 1311139 & 1311139 \\
wastewater05m1                      &   \bf{10.2} &  -1    &    -1  &    \bf{5519} &      -1 &      -1 \\
wastewater11m1                      & 3600.0 &    6.0 &    \bf{5.8} & 2913989 &    \bf{4315} &    \bf{4315} \\
wastewater13m1                      & 3600.0 &   11.0 &   \bf{10.9} & 1137029 &    \bf{1897} &    \bf{1897} \\
water4                              & 1515.8 & 1422.8 & \bf{1386.4} &  614339 &  572725 &  \bf{561704} \\
watercontamination0202r             & 3600.0 & 3600.0 &   \bf{33.8} &       \bf{5} &       \bf{5} &      13 \\
waterno1\_24                         &   \bf{20.8} & 3600.0 & 3600.0 &       \bf{1} &   19705 &   38321 \\
waterno2\_04                         &  \bf{631.3} & 3600.0 & 3600.0 &   \bf{61111} &  170213 &   79922 \\
watersym2                           & 3600.0 & 3600.0 & \bf{2607.0} & 1282308 & \bf{1082717} & 1161616 \\
\bottomrule
\end{tabular*}
\end{table}

\end{appendices}

\end{document}